 \documentclass[opre,nonblindrev]{informs3} 

\OneAndAHalfSpacedXI

\usepackage{subcaption}
 \usepackage{algorithm}
 \usepackage{algpseudocode}
\usepackage{mathrsfs}
\usepackage{amsmath, bm}
\usepackage{graphicx}
\usepackage{placeins}
 \usepackage{enumerate}
 \usepackage{fancyhdr}
\usepackage{url}
\definecolor{DarkBlue}{rgb}{0,0.08,0.45}
\usepackage[backref = false, bookmarks, colorlinks = true, plainpages = false, citecolor = DarkBlue , urlcolor = DarkBlue, filecolor = DarkBlue, linkcolor = DarkBlue]{hyperref}
\usepackage{newunicodechar}
\usepackage[utf8]{inputenc}

\usepackage{natbib}
\bibpunct[, ]{(}{)}{,}{a}{}{,}%
\def\bibfont{\small}%

\usepackage{pgfplots}
\usepgfplotslibrary{fillbetween}
\pgfplotsset{compat=1.13}
\usepackage{algorithm}
\usepackage{algpseudocode}
\usepackage{bbm}

\usepackage[color=yellow, open=true]{pdfcomment}

\usepackage{booktabs}
\usepackage{multirow}

\usepackage{tikz}
\usetikzlibrary{arrows.meta,positioning,shapes.geometric}

\newcommand{\changesh}[1]{#1}
\definecolor{orange}{rgb}{0,0,0}
\TheoremsNumberedThrough     
\ECRepeatTheorems

\EquationsNumberedThrough    

\begin{document}
\newcommand{\abs}[1]{\left|  #1 \right| }
\newcommand{\brak}[1]{\left(#1\right)}    
\newcommand{\crl}[1]{\left\{#1\right\}}   
\newcommand{\edg}[1]{\left[#1\right]}     
\newcommand{\norm}[1]{\|#1\|}
\newcommand{\floor}[1]{\lfloor #1 \rfloor}

\newcommand{\cA}{{\mathcal A}}
\newcommand{\cB}{{\mathcal B}}
\newcommand{\cD}{{\mathcal D}}
\newcommand{\cF}{{\mathcal F}}
\newcommand{\cG}{{\mathcal G}}
\newcommand{\cH}{{\mathcal H}}
\newcommand{\cK}{{\mathcal K}}
\newcommand{\cL}{{\mathcal L}}
\newcommand{\cM}{{\mathcal M}}
\newcommand{\cR}{{\mathcal R}}
\newcommand{\cS}{{\mathcal S}}
\newcommand{\cT}{{\mathcal T}}
\newcommand{\cX}{{\mathcal X}}
\newcommand{\cP}{{\mathcal P}}
\newcommand{\cV}{{\mathcal V}}

\newcommand{\mA}{{\mathbb A}}
\newcommand{\mV}{{\mathbb V}}
\newcommand{\mC}{{\mathbb C}}
\newcommand{\mR}{{\mathbb R}}
\newcommand{\mE}{{\mathbb E}}
\newcommand{\mw}{{\mathbb w}}
\newcommand{\mT}{{\mathbb T}}

\newcommand{\bb}{{\mathbf b}}
\newcommand{\bd}{{\mathbf d}}
\newcommand{\by}{{\mathbf y}}
\newcommand{\bp}{{\mathbf p}}
\newcommand{\bc}{{\mathbf c}}
\newcommand{\bg}{{\mathbf g}}
\newcommand{\bl}{{\mathbf l}}
\newcommand{\bbf}{{\mathbf f}}
\newcommand{\bq}{{\mathbf q}}

\newcommand{\bx}{{\mathbf x}}
\newcommand{\bA}{{\mathbf A}}
\newcommand{\bB}{{\mathbf B}}
\newcommand{\bC}{{\mathbf C}}
\newcommand{\bD}{{\mathbf D}}
\newcommand{\bG}{{\mathbf G}}
\newcommand{\bL}{{\mathbf L}}
\newcommand{\bS}{{\mathbf S}}
\newcommand{\bQ}{{\mathbf Q}}
\newcommand{\bU}{{\mathbf U}}
\newcommand{\bV}{{\mathbf V}}
\newcommand{\bX}{{\mathbf X}}
\newcommand{\bZ}{{\mathbf Z}}
\newcommand{\bF}{{\mathbf F}}

\newcommand{\bmu}{{\boldsymbol \mu}}
\newcommand{\bw}{{\boldsymbol w}}
\newcommand{\balpha}{{\boldsymbol \alpha}}
\newcommand{\blambda}{{\boldsymbol \lambda}}

\newcommand{\btheta}{\boldsymbol{\theta}}
\newcommand{\bsigma}{\boldsymbol{\sigma}}
\newcommand{\bnu}{\boldsymbol{\nu}}
\newcommand{\bSigma}{\boldsymbol{\Sigma}}
\newcommand{\bgamma}{\boldsymbol{\gamma}}
\newcommand{\bs}{\boldsymbol{s}}
\newcommand{\bz}{\boldsymbol{z}}

\newcommand{\F}{\mathbb{F}}
\newcommand{\p}{\mathbb{P}}
\newcommand{\Q}{\mathbb{Q}}

\newcommand{\R}{\mathbb{R}}
\newcommand{\q}{\mathbb{Q}}

\newcommand{\tr}{{\rm tr}}

\newcommand{\id}{{\mathbbm 1}}

\newcommand{\expect}{\mathbb{E}}



\RUNAUTHOR{Li, Peng, Shao, and Teo}

 \RUNTITLE{Extreme-Case Distorted Utility under Moment Ambiguity}

\TITLE{Extreme-Case Distorted Utility under Moment Ambiguity}

\ARTICLEAUTHORS{%
\AUTHOR{Zehao Li}
\AFF{Guanghua School of Management, Peking University, Beijing, 100871, China. zehaoli@stu.pku.edu.cn}

\AUTHOR{Yijie Peng}
\AFF{School of  Management and Engineering, Nanjing University, Nanjing, 210008, China. pengyijie@pku.edu.cn}

\AUTHOR{Hui Shao}
\AFF{International Business School, Zhejiang University, Haining, 314400, China. mathshao@gmail.com}

\AUTHOR{Chung-Piaw Teo}
\AFF{Institute of Operations Research and Analytics,
National University of Singapore, 117602, Singapore. bizteocp@nus.edu.sg }

} 
\ABSTRACT{
Many operations decisions under distributional ambiguity, from pricing and inventory to capacity and contracting, evaluate an action through a tail-sensitive distorted utility of an uncertain payoff and hedge against the least favorable distribution consistent with a few known moments; the resulting worst-case evaluation is the inner problem of a moment-based distributionally robust decision. We study this inner problem, the extreme-case distorted utility under moment constraints, for \textcolor{orange}{ a locally Lipschitz utility} that may be nonsmooth and neither convex nor concave together with a general, possibly atomic, distortion. Recasting the problem in the quantile domain, we develop a unified method that yields exact first-order optimality conditions and closed-form extremal values and distributions for both the worst and best cases, drawing on nonsmooth variational analysis. A central step treats the monotonicity constraint by isotonic projection onto the monotone cone, turning an abstract infinite-dimensional restriction into an inexpensive inner solve that scales linearly in the \textcolor{orange}{ discretization. 
The} method recovers and extends classical moment bounds through \changesh{three} examples: a range value-at-risk extension of the Scarf bound, GlueVaR distortions with a reward--penalty utility, and a capped incentive contract under conditional value-at-risk. As the inner oracle of a robust min-max decision, the characterization embeds directly in outer robust optimization, illustrated on a real capacity-provisioning problem for generative artificial intelligence inference where accounting for moment ambiguity lowers required capacity while preserving service compliance.
}%




\maketitle

\section{Introduction}

Many operational decisions optimize not the raw uncertain outcome but a nonlinear transformation of it, such as a service-level penalty, an option-like payoff, or a capped cost, aggregated so as to emphasize adverse tail events, while the distribution of the outcome is known only through a few moments. This paper studies the extreme-case distorted utility (DU), the worst-case or best-case value of such a transformed payoff over all distributions consistent with the given moments, a functional that generalizes classical risk and preference models under moment-based distributional ambiguity. Given a random variable $X$ with generalized quantile function $F_X^{-1}$, we define the DU as the integral of an increasing utility $U:\mathbb{R}\to\mathbb{R}$ with respect to  a nondecreasing distortion function $\varphi:\left [0,1\right ]\to[0,1]$. \changesh{We consider} the following optimization problem:
\begin{equation}\label{eq0}
    \operatorname*{\inf \ \text{or}\ \sup}_{X\in \mathcal{P}} \ \ 
    M_{U,\varphi}(X):=\int_{0}^{1}U\!\left(F_X^{-1}(u)\right)\mathrm{d}\varphi(u).
\end{equation}
Here, $\mathrm{d}\varphi$ denotes the Stieltjes measure induced by $\varphi$ and $\mathcal{P}$ is an ambiguity
set of square-integrable random variables, equivalently their laws, with moment information. The measure $\mathrm{d}\varphi$
reweights the utility values across quantile levels, thus allowing tail-sensitive or otherwise
nonlinear emphasis on different parts of the outcome distribution. Problem \eqref{eq0} subsumes several
canonical models: when $U(x)=x$, it reduces to the extreme-case distortion risk measure (DRM)
problem \citep{li2018closed, shao2023distortion, shao2024extreme}; when $\varphi(u)=u$, it reduces to optimizing the expected utility $\mathbb{E}[U(X)]$ over $\mathcal{P}$, i.e., a
generalized moment problem that also arises as the inner \changesh{step of} moment-based
distributionally robust optimization (DRO) \citep{smith1995generalized, bertsimas2011theory, kuhn2025distributionally}.

Our interest in problem~\eqref{eq0} is rooted in robust decision making: when an operational action such as a price, an order quantity, a reserved capacity, or a contract term shapes the payoff $X$, problem~\eqref{eq0} is the inner problem of the DRO decision that evaluates the action against the worst case over the distributions consistent with the known moments. Our contribution is to solve this inner problem in closed form for a general pair of utility and distortion.

Quantile-based risk measures such as value-at-risk (VaR) and conditional value-at-risk (CVaR) are canonical members of the broad family of tail-sensitive criteria represented as DRMs \citep{wang2000class}, which aggregate loss quantiles using a nondecreasing probability reweighting \citep{artzner1999coherent, dhaene2012remarks,iancu2015tight}.
This quantile-weighting perspective also has decision-theoretic roots in Choquet expected utility \citep{schmeidler1989, tversky1992cumulative}, a connection we make precise below.
 This representation subsumes many practically used specifications (e.g., spectral risk measure, range value-at-risk (RVaR), glue value-at-risk (GlueVaR), and parametric power-type distortions) while preserving a transparent interpretation in
terms of “which parts of the distribution matter.” In applications, however, the loss
distribution is rarely known with precision: one typically has finite data and only a small number of
reliably estimable summaries.  This motivates {extreme-case} evaluations over ambiguity sets that encode partial information, with moment constraints being particularly attractive due to their statistical
tractability and their role as core primitives in DRO \citep{popescu2007robust,zymler2013worst,chen2021discrete, pichler2022quantitative,shao2023distortion, shao2024extreme}.

Beyond unifying risk measures and rank-dependent preferences, the DU captures an economic distinction that classical \changesh{expected-utility and risk-measure models} cannot express. Under expected utility, $\mathbb{E}[U(X)]=\int_0^1 U(F_X^{-1}(u))\,\mathrm{d}u$ weights  quantile levels equally, and attitudes toward uncertainty are carried entirely by the curvature of $U$. The distortion $\varphi$ adds a second dimension: $U$ encodes preferences over outcomes, while $\varphi$ encodes preferences over the location of those outcomes within the distribution \citep{pesenti2024optimizing, cai2025distributionally}. This pairing is exactly the rank-dependent utility (RDU) of \cite{quiggin1982}\changesh{:} with a probability weighting $w$ inducing $\varphi(u)=1-w(1-u)$, the DU $M_{U,\varphi}(X)$ is the RDU of $X$, and Yaari's dual theory \citep{yaari1987} is the case $U=\mathrm{id}$\changesh{. The} extreme-case problem~\eqref{eq0} is \changesh{therefore} a DRO formulation of RDU: a decision maker who holds fixed RDU preferences but, lacking the full law, acts against the least-favorable distribution consistent with the moments. The separation matters in operations, where decision makers \changesh{overweight} specific regions of the distribution, such as \changesh{inventory stockouts}, service-level violations in capacity planning, \changesh{rare supply-chain disruptions}, or extreme losses in finance and insurance.

Distortion and ambiguity together yield a qualitatively different robust problem. Under DRO with expected utility, the adversary chooses a distribution within the ambiguity set to lower the average transformed payoff. Under DU, the adversary additionally allocates probability mass across quantiles that the distortion weights unequally, so ambiguity concerns not only the magnitude of uncertainty but also its location within the outcome distribution. The extreme-case distribution then depends jointly on the utility, the distortion function, and the available moments.

This interaction raises new economic questions. First, how does ambiguity affect a decision maker who is sensitive to a particular tail? Two agents with the same utility but different distortions can face entirely different least-favorable distributions. Second, how do moment constraints interact with tail emphasis? Under expected utility, variance bounds the overall dispersion of outcomes; under DU, the same variance has very different consequences depending on whether it sits in heavily weighted quantile regions. The economic value of reducing uncertainty thus depends not only on how much dispersion there is, but on where that dispersion falls in the distribution.

The extreme-case distributions themselves take a richer form. Classical moment problems typically yield discrete worst cases, readable as a small collection of stress scenarios; DU can produce discrete or continuous extremal distributions depending on how $U$ and $\varphi$ interact. Economically, this separates robustness against a finite set of identifiable adverse scenarios from robustness against diffuse distributional misspecification, and our characterization identifies when each regime arises.

The framework is especially suited to operational objectives that are naturally nonsmooth, such as threshold effects, service-level requirements, bonus and penalty clauses, inventory stockouts, option-like payoffs, or capped rewards, which generate piecewise-linear or otherwise nonsmooth utilities. In such settings,  the extreme-case distribution may strategically concentrate probability mass around operational thresholds, behavior that smooth expected-utility models cannot reproduce. Distortions with atoms arise just as naturally when attention focuses on specific quantiles, as in regulatory stress tests, service-level targets, or \changesh{VaR} constraints; such atoms induce an extreme-case {quantile targeting}, in which the least-favorable distribution concentrates its adverse effects at exactly the quantiles that carry explicit decision weight. Both features place the problem outside the convex-duality and smooth-optimization arguments used in existing theory.

\textbf{Literature review.} \changesh{Our work sits at the intersection of two streams of literature: the extreme-case evaluation of quantile-based risk measures under moment information, and DRO. The DU connects them: it reduces to a DRM when the utility $U$ is the identity and to an expected utility when the distortion $\varphi$ is the identity. The latter case is precisely the inner distributional problem of a moment-based DRO, which is why the two streams meet in our formulation.}

\changesh{The first stream concerns extreme-case DRMs, namely the worst- and best-case values of a quantile-based risk measure, together with the distributions that attain them, when the underlying law is known only through partial information such as a small number of moments. Foundational results on moment bounds and the associated extremal distributions appear in \cite{popescu2005semidefinite, popescu2007robust}. A substantial body of work then characterizes worst-case VaR and CVaR under moment-based ambiguity and related \textcolor{orange}{uncertainty descriptions \citep{ghaoui2003worst, zhu2009worst, natarajan2010tractable, chen2011tight, zymler2013worst,natarajan2018asymmetry}.} Moving beyond VaR and CVaR, \cite{cornilly2018upper} and \cite{li2018closed} derive closed-form worst-case values for broad classes of DRMs under moment constraints, covering strictly concave distortions and law-invariant coherent risk measures, respectively. \cite{shao2023distortion, shao2024extreme} unify and extend these results to general distortions through envelope-based constructions, and \cite{pesenti2024optimizing} establishes a convexification principle by which a nonconvex DRM can be optimized through its convex envelope, subject to verifiable structural conditions. A parallel line studies distortion-based robustness under alternative ambiguity models, notably Wasserstein-type uncertainty and \textcolor{orange}{related stability analyses \citep{mohajerin2018data,xie2021distributionally, liu2022inf, pichler2022quantitative, bernard2024robust},} as well as model aggregation grounded in stochastic dominance \citep{mao2025model}. Because the extreme-case DRM is the special case of the DU in which $U$ is the identity, these results are subsumed by our general characterizations.}

\changesh{The second stream is DRO, in which an extreme-case evaluation of the form $\sup_{X\in\mathcal{P}}\mathbb{E}\!\left [f(X)\right ]$ over an ambiguity set $\mathcal{P}$ arises as the inner problem of a min-max decision model. The classical Scarf bound is the canonical instance, giving a tight worst-case expected cost under mean--variance ambiguity \citep{scarf1958min, jagannathan1977minimax}. Subsequent work enriches the ambiguity description, for instance through higher-order moments and dependence \textcolor{orange}{structures \citep{delage2010distributionally,yu2022multistage, jiang2025dual}.} Closer to the DU, \cite{cai2025distributionally} show that a broad class of DRO problems with distortion-based objectives can be solved by replacing the distortion with its convex envelope, which yields tractable convex reformulations under additional structural conditions, and \cite{cai2024worst, cai2025conditional} evaluate extreme-case distorted utilities for specific parametric payoffs arising in insurance and asset management.}

\changesh{Taken together, the two streams obtain tractability by restricting either the utility or the distortion. The utility is the identity (as in the DRM bounds), concave (so that the objective can be convexified through the convex envelope of the distortion), or a fixed parametric payoff; and the distortion is, with few exceptions, required to be atomless. To the best of our knowledge, no prior work characterizes the extreme case for a general utility that may be neither convex nor smooth, together with a general and possibly atomic distortion, under moment information. We close this gap by providing exact, closed-form characterizations of both the worst and best cases. Because the extreme-case DU is exactly the inner problem of a moment-based DRO, these characterizations serve as reusable analytical building blocks for min-max decision models whose objectives combine probability distortion with outcome transformation, including robust pricing and revenue management, robust capacity provisioning, and robust portfolio and reinsurance design. Table~\ref{tab:related_special_cases} summarizes the related work and positions our study within this literature.} 

\begin{table}[htbp]
\centering
\footnotesize
\setlength{\tabcolsep}{2pt}
\renewcommand{\arraystretch}{1.00}
\begin{tabular}{lcccc}
\toprule
\textbf{Reference} 
& \textbf{Distortion function} 
& \textbf{Utility } 
& \textbf{Worst-case} 
& \textbf{Best-case} \\
\midrule
\cite{scarf1958min}, \cite{jagannathan1977minimax}
& $\times$ 
& piecewise linear 
& \checkmark 
& $\times$ \\

\cite{ghaoui2003worst}
&   VaR
& $\times$   
& \checkmark
& $\times$ \\

\cite{popescu2005semidefinite}
& $\times$
& general payoff 
& \checkmark
& \checkmark \\

\cite{natarajan2010tractable}
& $\times$
& piecewise-linear concave
& \checkmark 
& $\times$ \\

\cite{chen2011tight}
& CVaR
& $\times$
& \checkmark 
& $\times$ \\

\cite{li2018worst}
& RVaR
& $\times$
& \checkmark 
& $\times$ \\

\cite{li2018closed}
& concave
& $\times$
& \checkmark 
& $\times$ \\

\cite{cornilly2018upper}
& strictly concave
& $\times$
& \checkmark 
& $\times$  \\

\cite{liu2022inf}
& tail risk measures
& $\times$
& \checkmark 
& $\times$  \\

\cite{bernard2024robust}
& absolutely continuous
& $\times$
& \checkmark 
& \checkmark  \\

\cite{cai2024worst}
& general
& $\left (x-a\right )_+$, $x\wedge m$
& \checkmark
& $\times$ \\

\cite{shao2023distortion, shao2024extreme}
& \multirow{2}{*}{general}
& \multirow{2}{*}{$\times$}
& \multirow{2}{*}{\checkmark}
& \multirow{2}{*}{\checkmark} \\
\cite{pesenti2024optimizing}
& 
& 
& 
& \\ 

\cite{liucoverage}
& general
& $\times$
& \checkmark 
& \checkmark  \\

\cite{jiang2025dual}
& $\times$ 
& $\left (\sum_i x_i-a\right )_+$
& \checkmark 
& $\times$ \\

\cite{cai2025distributionally}
& ``insensitive''
& concave
& \checkmark 
& $\times$ \\

\cite{cai2025conditional} 
& CVaR
& $x+\kappa\left (x-a\right )_+$
& \checkmark
& $\times$ \\

\hline

This paper 
& general  
& locally Lipschitz
& \checkmark
& \checkmark \\



\bottomrule
\end{tabular}
\caption{Positioning of this paper within the extreme-case literature. Each row lists the distortion, utility, and whether worst and best cases are treated. A check mark ($\checkmark$) indicates that the component is addressed, and a cross ($\times$) denotes an identity utility or distortion. Here $\left (x\right )_+:=\max\!\left \{x,0\right \}$ and $x\wedge m:=\min\!\left \{x,m\right \}$. An ``insensitive'' distortion \changesh{vanishes on $[0,\epsilon]$, the condition under which} the convex-envelope reformulation is exact.}
\label{tab:related_special_cases}
\end{table}

\changesh{A related but distinct stream studies preference ambiguity, optimizing over a set of plausible utility or distortion functions rather than over the loss distribution \citep{wang2023preference, guo2024utility, wang2025preference}. We take the complementary view, fixing the pair $(U,\varphi)$ and characterizing the extreme cases induced by moment ambiguity in the distribution. The two formulations combine when distributional and preference ambiguity are modeled jointly.}

\textbf{Contributions.} This paper makes three contributions.

(1) \changesh{Exact, closed-form characterization of extreme-case DU. For a general pair $(U,\varphi)$, with a utility that need not be concave, convex, or smooth and a distortion that may carry atoms, we derive exact first-order optimality conditions for both the worst- and best-case problems and obtain the extremal values, together with the distributions that attain them, in closed form. Methodologically, we work directly with the original nonconvex and nonsmooth problem rather than a convexified or relaxed surrogate \citep{shao2023distortion, shao2024extreme, pesenti2024optimizing, cai2025distributionally}: a quantile-domain variational analysis encodes the monotonicity of the quantile function as a normal-cone condition and resolves it constructively by isotonic projection onto the monotone cone.}

(2) \changesh{A unified and tractable solution framework. The optimality system reduces to a single algorithm that applies across a broad class of objectives, with a closed-form inner solve whose cost is linear in the \textcolor{orange}{discretization.} We further establish existence and uniqueness of the extreme-case optimizer under mild regularity conditions.} 
\changesh{The framework recovers and extends classical moment bounds, which we illustrate through three representative cases: an RVaR extension of the Scarf bound, GlueVaR with a reward--penalty utility, and a capped incentive contract under CVaR provisioning.}

(3) \changesh{Structural and operational insights for robust decisions. The characterization reveals how the shape of the extremal distribution is governed by the interplay of utility and distortion, ranging from low-dimensional discrete laws for piecewise-linear payoffs to continuous distributions for smooth utilities, and it identifies regimes in which the extremum is not attained. It further traces how the extreme-case value responds to the decision and cost parameters, yielding transparent comparative statics. Because the extreme-case DU is exactly the inner problem of a moment-based DRO, these characterizations serve as reusable closed-form inner oracles that can be evaluated repeatedly within the outer min-max loop of robust decisions such as pricing, inventory, reinsurance, and capacity provisioning. We develop a real-data application to robust capacity provisioning for generative artificial intelligence inference in Section~\ref{sec:genai-ladder} and quantify the value of accounting for ambiguity.}


The rest of the paper is organized as follows. Section \ref{sec2} introduces the DU formulation. Section \ref{sec3} establishes the optimality theory for extreme-case DU. Section \ref{sec4} develops the solution algorithm, closed-form examples, and the real-data robust decision application. Section \ref{sec5} concludes. Proofs and supplementary results are in the online appendices.

\section{Problem Statement}\label{sec2}
This section formalizes the DU and associated extreme-case problems
under the first two moment constraints. 
Let $(\Omega,\mathcal F,\mathbb P)$ be a probability space and $X:\Omega\to\mathbb R$ be a random variable with cumulative distribution function (CDF) $F_X$, i.e., $F_X(t):=\mathbb{P}(X\le t)$ for $t\in\mathbb{R}$. The generalized inverse (quantile) function of $X$ is defined by $F_X^{-1}(u):=\inf\{t\in\mathbb R: F_X(t)\ge u\}$ for $u\in(0,1)$. Throughout this paper, we assume that all the random variables are square-integrable, i.e., $X\in L^2\left (\Omega, \mathcal{F}, \mathbb{P}\right )$. Let $L^2(0,1)$ be a real Hilbert space with the inner product
$\left \langle q_1,q_2\right \rangle:=\int_0^1 q_1(u)q_2(u)\,\mathrm{d}u$ and norm $\left \|q\right \|_2:=\left (\int_0^1 q^2(u)\,\mathrm{d}u\right )^{1/2}$  for any $q, q_1, q_2\in L^2(0,1)$.  

\begin{definition}[DU]\label{def:DU}
\changesh{Let $U:\mathbb R\to\mathbb R$ be an increasing function (not necessarily concave, convex, or smooth) and  $\varphi$ be a distortion function belonging to the family}
\begin{displaymath}
\mathcal{D}:=\left \{\varphi: [0,1]\to [0,1]\mid \varphi(0)=0,\, \varphi(1)=1,\, \text{$\varphi$ is nondecreasing and right-continuous} \right \}. 
\end{displaymath}
\changesh{For any random variable $X$ for which the following integral is finite, the DU induced by the pair $(U,\varphi)$ is}
\begin{displaymath}\label{eq:DU}
M_{U, \varphi}(X)
:=\int_{0}^{1} U\!\left (F_X^{-1}(u)\right ) \mathrm d\varphi(u).
\end{displaymath}
\changesh{Here $\mathrm{d}\varphi$ is the Stieltjes measure induced by $\varphi$, and $F_X^{-1}$ denotes the quantile function of $X$.}
\end{definition}

In particular, when either $U$ or $\varphi$ is the identity, we obtain well-known forms of utilities or risk measures. When the distortion is the identity, $\varphi(u)=u$, the DU is the expected utility $M_{U,\mathrm{id}}(X)=\mathbb{E}[U(X)]$ \citep{neumann1944, savage1954}, with $U$ valuing the outcomes; when instead the utility is the identity, the DU is the conventional DRM $M_{\mathrm{id},\varphi}(X)=\int_0^1 F_X^{-1}(u)\,\mathrm{d}\varphi(u)$ \citep{yaari1987, wang2000class}, in which $\varphi$ reweights the quantiles, and specializing $\varphi$ recovers standard risk measures such as VaR and CVaR.

\begin{definition}[Moment Constraints]\label{def2}
Given a mean $\mu\in\mathbb R$ and a standard deviation $\sigma\ge 0$, let
\begin{displaymath}
\mathcal P(\mu,\sigma^2)
:=\left\{ X\in L^2(\Omega, \mathcal{F}, \mathbb{P})\mid  \mathbb{E}\!\left [X\right ]=\mu, \mathbb{E}\!\left [X^2\right ]=\mu^2+\sigma^2
    \right\}
\end{displaymath} be the set of random variables that match the specified first two 
moments. Accordingly, the worst-case and best-case DUs can be defined as 
\begin{equation}\label{extreme_case_problem}
\sup_{X\in\mathcal{P}(\mu,\sigma^2)} M_{U,\varphi}(X)\quad\text{and}\quad \inf_{X\in\mathcal{P}(\mu,\sigma^2)} M_{U,\varphi}(X), 
\end{equation}respectively. 
\end{definition}

We are concerned with the extreme-case DU defined in problem \eqref{extreme_case_problem}. 
Owing to the generality of the DU, problem~\eqref{extreme_case_problem} includes several well-established problems in the literature as special cases.  


Although the extremal problem \eqref{extreme_case_problem} is posed on the space of random variables, it is analytically convenient to recast it in the quantile domain. We define the {monotone cone} as
\begin{displaymath}\mathsf K:=\left \{q\in L^2(0,1): q \text{ admits an almost everywhere (a.e.) nondecreasing representative on }(0,1)\right \}.\end{displaymath} A function $q:\left (0,1\right )\to\mathbb R$ is nondecreasing a.e.\ if it admits a representative, still denoted by $q$, that is nondecreasing on $(0,1)$ outside a Lebesgue null set. 
We use this representative to define the nonnegative Stieltjes measure $\mathrm{d}q$.  Using the standard quantile-integral representation in \cite{dhaene2012remarks}, the problem can be equivalently written as
\begin{equation}\label{p1}
\begin{aligned}
& \underset{}{\inf_q\ \text{or}\ \sup_q}
& & \Psi(q):=\int_0^1 U (q(u))\, \mathrm{d}\varphi (u) \\
& \text{subject to}
& & \int_0^1 q(u)\,\mathrm{d}u=\mu,\,\, \int_0^1 q^2(u)\,\mathrm{d}u=\mu^2+\sigma^2,\,\,   q\in\mathsf{K}. 
\end{aligned}
\end{equation}
Here, $q$ denotes the quantile function (the decision variable), and $\mathsf{K}$ enforces the monotonicity of the admissible quantile functions. The feasible set is nonconvex, owing to the second-moment equality constraint. Problem \eqref{p1} is a nonsmooth functional optimization problem, to which we apply the Clarke subdifferential developed for locally Lipschitz functions. 
Next, we introduce definitions from Clarke's nonsmooth analysis, which will be used in our later developments.

\begin{definition}[\cite{clarke1990optimization}]
Let $\mathcal X$ be a real Banach space.
For a locally Lipschitz function $f:\mathcal X\to\mathbb R$, the {Clarke directional derivative} of $f$ at $x\in \mathcal X$ in direction $h\in \mathcal X$ is defined as
\[
f^\circ(x;h):=\limsup_{\substack{y\to x,  t\to 0^+}} \frac{f(y+t\, h)-f(y)}{t}.
\]
The {Clarke subdifferential} of $f$ at $x$, also called the Clarke generalized gradient, is defined as
\[
\partial f(x):=\left \{v\in \mathcal X^*: \left \langle v,h\right \rangle \le f^\circ(x;h),\, \forall h\in \mathcal X\right \},
\]
where $\mathcal X^*$ denotes the dual space of $\mathcal X$. When $\mathcal X=L^2(0,1)$, we identify $\mathcal X^*$ with $\mathcal X$ by the Riesz representation theorem \citep{rudin1976principles}.
\end{definition}

Intuitively, the Clarke directional derivative generalizes one-sided derivatives, and $\partial f$ serves as the gradient in nonsmooth first-order conditions.

\section{Main Results}\label{sec3}
This section presents our main theoretical results. 
%
To solve the optimization problem \eqref{p1}, we impose  the following mild assumptions on the pair of utility and distortion functions in the DU: 
\begin{assumption}\label{asp1}
$\varphi:[0,1]\to [0,1]$ is absolutely continuous with $\varphi' \in L^\infty(0,1)$.
\end{assumption}
\begin{assumption}\label{asp2}
$U:\mathbb R\to\mathbb R$ \textcolor{orange}{is locally} Lipschitz, i.e., for each $x_0\in\mathbb{R}$ there exist constants $\delta_0>0$ and $M>0$ such that 
$\left |x-x_0\right |<\delta_0\implies \left | U(x)-U(x_0)\right |\le M\! \left |x-x_0\right |. $
\end{assumption}
\begin{assumption}\label{asp3}
There exist $a,b\ge 0$ such that
$\sup_{\zeta\in\partial U(x)} |\zeta| \le a+b\left |x\right |,\,\, \forall\,x\in\mathbb R,$
where $\partial  U(x)$ denotes the Clarke subdifferential of $U$ at $x$. 
\end{assumption}

Assumption \ref{asp1} ensures the global Lipschitz continuity of $\varphi$. Note that this condition is subsequently relaxed to accommodate general distortion functions that are not necessarily absolutely continuous. Assumption \ref{asp2} establishes the existence of the Clarke subdifferential for the functional formulated in problem \eqref{p1}. Finally, Assumption~\ref{asp3} restricts the utility function to at most quadratic growth; this guarantees the well-posedness of both the functional and its Clarke directional derivatives in $L^2$, while preserving the model's applicability to nonsmooth and nonconcave utilities.

We work with a nonsmooth optimality framework that explicitly encodes the monotonicity constraint $q\in\mathsf K$. 
Instead of treating $\mathsf K$ as an external restriction, we incorporate it into the objective by adding the indicator functional $ \mathcal I_{\mathsf K}$ \changesh{defined by}
$$\changesh{\mathcal I_{\mathsf K}(q):=\begin{cases} 0, & q\in\mathsf K,\\[3pt] +\infty, & \text{otherwise.} \end{cases}}$$ 
With this reformulation, the effect of the constraint is captured entirely through the generalized subdifferential of $ \mathcal I_{\mathsf K}$.
In particular, the subdifferential of the indicator is characterized by the {normal cone} to $\mathsf K$: $\partial \mathcal I_{\mathsf K}(q)=N_\mathsf K (q) $. 

\begin{definition}[Polar and Normal Cone]\label{def1} 
The polar of a monotone cone $\mathsf K$ is defined as 
$\mathsf K^\circ:=\left \{z\in L^2:  \langle z,h\rangle\le 0,\, \forall h\in\mathsf K\right \},$
and the normal cone at $q\in\mathsf K$ is defined as
$N_{\mathsf K}(q):=\left \{z\in L^2: \left \langle z, h-q\right \rangle\le 0,\, \forall h\in \mathsf K\right \}.$
\end{definition}

The polar cone collects the functions whose inner product with every monotone function is nonpositive, and the normal cone at a monotone $q$ collects those whose inner product with $h-q$ is nonpositive for every $h\in\mathsf K$, a form of orthogonality. By  definition of $\mathsf K$, the polar and normal cones at $q\in\mathsf K$ admit the following sharper characterization.

\begin{proposition}[Normal {cone} of $\mathsf K$]\label{lemma:normal_cone_K} 
For any $q\in L^2(0,1)$, define its primitive function as 
$\mathcal{J}\!q(t):=\int_0^t q(s)\,\mathrm{d}s$ for $t\in[0,1]$. 
Then the following conclusions hold.  
\begin{enumerate}
\item[(a)]  The polar of the monotone cone $\mathsf{K}$ is
$\mathsf K^\circ=\left \{z\in L^2(0,1): \mathcal{J}\!z(t) \ge 0,\, \forall t\in[0,1],\,\mathcal{J}\!z(1)=0\right \}.$

\item[(b)]  For any $q\in\mathsf K$, the normal cone at $q$ is 
\begin{equation}\label{optimality condition2}
N_{\mathsf K}(q)=\mathsf K^\circ\cap\left \{q\right \}^\perp
=\left \{z\in L^2(0,1): \mathcal{J}\!z(t) \ge 0,\, \forall t\in[0,1],\, \mathcal{J}\!z(1)=0,\, \int_{0}^1  \mathcal{J}\!z(t)\,\mathrm{d}q(t)=0\right \},
\end{equation}
where $\left \{q\right \}^\perp:=\left \{z\in L^2(0,1): \left \langle z,q\right \rangle=0\right \}$. 
\end{enumerate}
\end{proposition}

Proposition~\ref{lemma:normal_cone_K} shows that the monotonicity constraint can be handled
through a simple complementarity structure. The primitive $\mathcal{J}\!z$ acts as a multiplier that is nonnegative everywhere. The complementarity condition $\int_0^1 \mathcal{J}\!z(t)\,\mathrm{d}q(t)=0$ indicates that the primitive of any normal-cone element vanishes \changesh{where $q$ strictly increases and is} positive only where $q$ is constant.

Then, we formalize the moment constraints as an equality set $\mathsf M$ for a given mean $\mu$ and standard deviation $\sigma$. Corresponding to Definition \ref{def2}, we define the closed constraint set $\mathsf M$ as
$$\mathsf M:=\left \{q\in L^2(0,1): \left \langle 1,q\right \rangle=\mu,\left  \langle q,q\right \rangle=\mu^2+\sigma^2\right \}.$$
At an optimal point $q^*\in\mathsf M$, the moment constraints are equalities and thus $\mathsf M$ is generally neither a cone nor convex, unlike the monotonicity constraint set $\mathsf K$ considered earlier. For this reason, instead of using polar-cone geometry as we did for $\mathsf K$, we describe first-order feasible perturbations of $q^*$ through the Bouligand (contingent) tangent cone $T_{\mathsf M}(q^*)$, defined as $$T_{\mathsf M}(q^*)
:=\left \{H\in L^2(0,1): \text{ there exists } t_n\downarrow 0\ \text{and } q_n\in \mathsf M,  \text{ s.t.}\, 
\left \|q_n-\left (q^*+t_n H\right )\right \|_{2}=o(t_n)\right \}.$$ 

Under the linear independence constraint qualification (LICQ), namely that $1$ and $2\,q^*$ are linearly independent in $L^2$, the equality structure makes the tangent cone the intersection of the kernels of the linearized moment maps, and the normal cone the span of their gradients.

\begin{proposition}[Normal cones of $\mathsf M$]\label{lem:tangent-normal-M}
Let $q^*\in \mathsf M$ be the optimal point and assume LICQ at $q^*$, or equivalently, $\sigma^{2}>0$.
Then the  contingent tangent cone can be represented as  $T_{\mathsf M}(q^*)=\left \{H\in L^2: \left \langle 1,H\right \rangle= \left \langle 2\, q^*,H\right \rangle=0\right \}$, and the normal cone at $q^*$ is $N_{\mathsf M}(q^*)=\mathrm{span}\!\left  \{1,2\, q^*\right \}$.
\end{proposition}

Proposition~\ref{lem:tangent-normal-M} expresses the local geometry of $\mathsf{M}$ through its tangent and normal cones, providing a key link between the constraints and the optimality conditions. The following remark records a degenerate case that we rule out in the subsequent analysis.

\begin{remark}[Degenerate variance]\label{rmk:degenerate}
If $\sigma^{2}=0$, then $q^*$ is constant and $1,2q^*$ are linearly dependent (LICQ fails).
In this case, $\mathsf M=\left \{q^*\right \}$, and we obtain $T_{\mathsf M}(q^*)=\left \{0\right \}$ and $N_{\mathsf M}(q^*)=L^{2}$; {the conclusion
in Proposition~\ref{lem:tangent-normal-M} does not hold.}
\end{remark}

\subsection{Optimality conditions under absolutely continuous distortions} 

In this subsection, we derive the optimality conditions for problem \eqref{p1}. The following lemma verifies that the objective functional meets the basic requirement of nonsmooth analysis, namely local Lipschitz continuity, which is the prerequisite for applying the Clarke subdifferential in the derivation that follows.

\begin{lemma}[Local Lipschitzness of $\Psi$]\label{lem:Psi-Lip}
Under Assumptions \ref{asp1}, \ref{asp2}, and \ref{asp3}, the functional $\Psi:L^2(0,1)\to\mathbb R$ in  problem \eqref{p1} is locally Lipschitz: for each $q_\circ\in L^2$ and $r>0$,  there exists
$L_{q_\circ,r}:=\left \|\varphi^\prime \right \|_{\infty}\left (a+b\left (2r+2\left \|q_\circ\right \|_{2}\right )\right )$
such that for all $q_1,q_2\in L^2$ satisfying $\left \|q_i-q_\circ\right \|_{2}\le r$, we have
$\left |\Psi(q_1)-\Psi(q_2)\right | \le L_{q_\circ,r} \left \|q_1-q_2\right \|_{2}$.
\end{lemma}

For an absolutely continuous distortion function, $\mathrm d\varphi(u)=\varphi'(u)\,\mathrm du$ for some density $\varphi'$, with no atomic or singular part. Combining the two normal-cone results above, the next theorem gives the optimality condition for $q^*$: monotonicity appears through \eqref{optimality condition2} with $z^*\in N_{\mathsf K}(q^*)$, and the moment constraints through $N_{\mathsf M}(q^*)$.

\begin{theorem}[Optimality condition for minimization]\label{main_theorem1}
Suppose Assumptions \ref{asp1}, \ref{asp2}, and \ref{asp3} and LICQ at $q^*\in\mathsf M\cap\mathsf K$ hold. 
If $q^*$ is a local minimizer of $\Psi$ over $\mathsf M\cap\mathsf K$, then the following conclusions hold. All subgradients and normal-cone elements are understood as elements of $L^2$.

(a) There exists a pair of multipliers
$(\lambda_1,\lambda_2)$ such that
$0 \in \partial \Psi(q^*) + \lambda_1 + 2\,\lambda_2\,q^* + z^*$, where $z^*\in N_{\mathsf K}(q^*)$ is characterized by Equation \eqref{optimality condition2}.

(b) Moreover, there exists a measurable {selector} $\xi(\cdot)\in\partial U(q^*(\cdot))$ and $z^*\in N_{\mathsf K}(q^*)$ with
\begin{equation}\label{optimality condition}
0 \in \varphi^\prime\,\xi + \lambda_1 +2\, \lambda_2\,q^* + z^*.
\end{equation}

(c) On any open interval where $q^*$ is (a.e.) strictly increasing, the above conclusion (b) reduces to the pointwise condition: 
$-\lambda_1-2\lambda_2\,q^*(u) \in \varphi^\prime (u)\,\partial U\!\left (q^*(u)\right )$ for a.e. $u\in[0,1]$.
\end{theorem}

Theorem~\ref{main_theorem1} characterizes a local minimizer by a stationarity inclusion in which two Lagrange multipliers, for the mean and variance constraints, and a normal-cone element, for the monotonicity constraint, balance the subdifferential of the objective, with part (b) making $\partial\Psi(q^*)$ explicit. The key consequence is part (c): on intervals where $q^*$ is strictly increasing the normal-cone term vanishes, leaving the pointwise condition that drives the closed-form solutions of Section~\ref{sec4}.

The proof must handle a nonsmooth, only locally Lipschitz objective together with a feasible set formed by nonconvex moment equalities and an active monotonicity constraint, so classical pointwise Karush--Kuhn--Tucker (KKT) arguments do not apply. The key step characterizes $\partial\Psi$ in $L^2(0,1)$ by controlling the Clarke directional derivative and interchanging $\limsup$ with integration, which yields the almost-everywhere pointwise inclusion \citep{clarke1990optimization,clarke1998nonsmooth,rockafellar1998variational,mordukhovich2006variational}; this variational route differs from approaches based on H\"older's inequality, duality, or the convex envelope of nonconvex distortions \citep{popescu2007robust,gao2023distributionally,shao2024extreme,pesenti2024optimizing,cai2025distributionally,chung2026unbalanced,jiang2025dual}. 

Similarly, we have the following necessary condition for the maximization problem \eqref{p1}.
\begin{corollary}[Optimality condition for maximization]\label{corollary:OPT_for_max}
Suppose Assumptions \ref{asp1}, \ref{asp2}, \ref{asp3} and LICQ at $q^*\in\mathsf M\cap\mathsf K$ hold.
If $q^*$ is a local maximizer of $\Psi$ over $\mathsf M\cap\mathsf K$, then there exist multipliers
$(\lambda_1,\lambda_2)\in\mathbb R^2$, an element $z^*\in N_{\mathsf K}(q^*)$ and a measurable selector $\xi(\cdot)\in\partial U(q^*(\cdot))$, such that
\begin{equation*}
0\in-\varphi^\prime\,\xi+\lambda_1+2\lambda_2\,q^*+z^*.
\end{equation*}
On any open interval where $q^*$ is (a.e.) strictly increasing, the equation reduces to the pointwise condition
$\lambda_1+2\lambda_2\,q^*(u) \in \varphi^\prime(u)\,\partial U\!\left (q^*(u)\right )$ for a.e. $u\in[0,1]$. 
\end{corollary}


Similar to the minimization problem, {the maximization case differs only by a sign in the objective. The resulting optimality condition has the same structure: the negative subgradient term $-\varphi'\xi$ is balanced by a linear combination of constraint gradients and a normal-cone term for the monotonicity constraint. 
Importantly, the specific choice of the first two moment equalities in $\mathsf M$ is not essential for the derivation of the optimality condition above. 
What we use from $\mathsf M$ is only a local first-order description at $q^*$, namely that under \textup{LICQ} the tangent cone can be characterized by the linearized equalities and the associated normal cone is the span of the corresponding gradients. 
Hence, the same argument applies to other collections of moment equality constraints.

For this reason, Theorem~\ref{main_theorem1} extends to the first $n$ moment conditions; the statement is given in Appendix C. \textcolor{orange}{For the completeness of theory, we also establish the existence and uniqueness of optimizers for problem \eqref{p1} in Appendix A.8.}

\subsection{Optimality conditions under general distortion functions} 
This subsection explores the case where $\varphi$ is not absolutely continuous. By Lebesgue's decomposition theorem \citep{halmos1974measure}, the Stieltjes measure $\mathrm{d}\varphi$ decomposes as
$(\mathrm{d}\varphi)_{\rm ac}+\sum_{j=1}^J m_j\,\delta_{t_j}+(\mathrm{d}\varphi)_{\rm sc}$,
where  $t_j\in(0,1]$, $m_j>0$, $J$ is finite,
and $(\mathrm{d}\varphi)_{\rm sc}$ is singular-continuous. In most distortion models of practical interest, $(\mathrm{d}\varphi)_{\rm sc}=0$, and the jump part $\sum_{j=1}^J m_j\,\delta_{t_j}$ is the main focus of this subsection.  Let
$\mathrm{d}\varphi(u)=w(u)\,\mathrm{d}u+\sum_{j=1}^J m_j\,\delta_{t_j}(\mathrm{d}u)$ for 
$w\in L^\infty(0,1)$, $m_j> 0$, and $t_j\in(0,1]$.
For this distortion, define
\[
\Psi(q):=\int_0^1 w(u)\,U(q(u))\,\mathrm{d}u+\sum_{j=1}^J m_j\,U\!\left (q(t_j)\right )\,.
\]

The following theorem extends Theorem \ref{main_theorem1} to general distortions. The optimality condition now carries additional contributions at the jump points, while its remaining structure is unchanged, still balancing the Lagrange multipliers, a normal-cone element, and the subdifferential of the objective. \textcolor{orange}{
Before stating the theorem, \changesh{we} introduce the dominating measure
\(\mathrm d\nu:=\mathrm du+\sum_{j=1}^J\delta_{t_j}(\mathrm du)\) and the Hilbert space
\(\mathcal H:=L^2(\mathrm d\nu)\). The moment set \(\mathsf M\) is still defined through the ordinary
Lebesgue moments, namely \(\int_0^1 q(u)\,\mathrm du=\mu\) and
\(\int_0^1 q^2(u)\,\mathrm du=\mu^2+\sigma^2\). Let
\(\rho_\varphi:=\mathrm d\varphi/\mathrm d\nu\) and \(\rho_0:=\mathrm du/\mathrm d\nu\). Thus
\(\rho_\varphi=w\) \(\mathrm du\)-a.e. and \(\rho_\varphi(t_j)=m_j\), while
\(\rho_0=1\) \(\mathrm du\)-a.e. and \(\rho_0(t_j)=0\). All normal cones in this subsection are taken in
\(\mathcal H\), and for \(z\in\mathcal H\) we write
\(\mathcal J_\nu z(t):=\int_{[0,t]}z(s)\,\mathrm d\nu(s)\).  }

\begin{theorem}[Optimality condition for general DU]\label{main_thm2_general}
Under Assumptions \ref{asp2}, \ref{asp3}, and LICQ, 
if \(q^*\) is a local minimizer of \(\Psi\) over \(\mathsf K\cap\mathsf M\) with respect to the
\(\mathcal H\)-topology, then there exist multipliers
\((\lambda_1,\lambda_2)\in\mathbb R^2\), a measurable selector
\(\xi(u)\in\partial U(q^*(u))\) for \(\mathrm du\)-a.e. \(u\), pointwise selectors
\(\xi_j\in\partial U(q^*(t_j))\), and a normal term \(z^*\in N_{\mathsf K}(q^*)\) such that
\[
0=\rho_\varphi\widehat \xi+(\lambda_1+2\lambda_2 q^*)\rho_0+z^*
\  \text{in } \mathcal H,
\]
where \(\widehat\xi=\xi\) \(\mathrm du\)-a.e. and \(\widehat\xi(t_j)=\xi_j\).
Equivalently, the variational inequality holds for every \(V\in\mathsf K\):
\begin{equation}\label{eq:VI}
\int_0^1 \left (w(u)\xi(u)+\lambda_1+2\lambda_2 q^*(u)\right )\left (V(u)-q^*(u)\right ) \mathrm{d}u
+\sum_{j=1}^J m_j\xi_j\left (V(t_j)-q^*(t_j)\right )
\ge0.
\end{equation}
Moreover, the normal-cone element satisfies \(\mathcal J_\nu z^*(t)\ge0\) for \(t\in[0,1]\),
\(\mathcal J_\nu z^*(1)=0\), and
\(\int_0^1\mathcal J_\nu z^*(t)\,\mathrm dq^*(t)=0\), where the last integral is the Stieltjes integral.
\end{theorem}


\textcolor{orange}{
Theorem~\ref{main_thm2_general} covers practically relevant DRMs with jumps (e.g., GlueVaR) and yields a unified optimality system that simultaneously accounts for the continuous density part \(w\,\mathrm du\) and the discrete mass points \(\sum_{j=1}^J m_j\delta_{t_j}\). Equivalently, the variational inequality is a geometric first-order condition for the Lagrangian: the equality moment constraints are encoded by the multipliers \((\lambda_1,\lambda_2)\), while the monotonicity constraint is encoded by the normal-cone element \(z^*\in N_{\mathsf K}(q^*)\). Hence, after the moment constraints are dualized, no first-order descent direction remains within the monotone cone. Moreover, the primitive \(\mathcal J_\nu z^*\) satisfies a Stieltjes-type complementarity condition, implying that \(\mathcal J_\nu z^*\) can be positive only on plateau regions of \(q^*\); on strictly increasing parts, the normal-cone correction is inactive in the corresponding complementarity sense. }

\textcolor{orange}{
Compared with Theorem~\ref{main_theorem1}, the main technical difficulty is that the atomic objective involves point evaluations, which are not well defined, and hence not continuous, on \(L^2(\mathrm du)\). We therefore work in the dominating Hilbert space \(\mathcal H:=L^2(\mathrm d\nu)\), where both the atomic evaluations and the Lebesgue moment maps are continuous. The proof applies the Clarke multiplier rule, with Radon--Nikodym densities \(\rho_\varphi\) and \(\rho_0\) separating the distorted objective from moment constraints.  }

We close this subsection by placing the convex-envelope method of \cite{shao2024extreme} inside our framework. The next proposition shows that the quantile-projection stationarity is strictly more general: for the identity utility under mean--variance ambiguity it reproduces, in both value and optimizer, the closed-form envelope reduction of \cite{shao2024extreme}, so that the established convex-envelope theory is exactly the $U(x)=x$ special case of our analysis.

\begin{proposition}[The convex-envelope reduction as a special case]
\label{prop:recover-envelope-k2}
Let $\varphi^*$ be the convex envelope of the distortion function $\varphi$ in the sense of \cite{shao2024extreme}, i.e.,
\begin{equation*}
\varphi^*(u):=\sup\!\left \{g(u): g\le \varphi \ \text{on }[0,1],\ g \text{ is convex on }[0,1]\right \},\quad u\in[0,1].
\end{equation*}
Consider the special case $U(x)=x$ of our DU framework.
Assume further that $\varphi^*$ is nontrivial so that $(\varphi^*)'$ exists a.e.\ and
$(\varphi^*)'\not\equiv 1$.
Then the maximizer of $M_{U,\varphi^*}$ over $\mathsf K\cap\mathsf M$ is
\begin{equation}\label{eq:qstar-envelope}
q^*(u)
=\mu+\sigma\,\frac{(\varphi^*)'(u)-1}{\left \|(\varphi^*)'-1\right \|_2},
\quad u\in[0,1]\ \text{a.e.},
\end{equation}
and the optimal value is
\begin{equation}\label{eq:opt-envelope}
\sup_{q\in\mathsf K\cap\mathsf M}M_{U,\varphi^*}(q)
= M_{U,\varphi^*}(q^*)
= \mu+\sigma\left \|(\varphi^*)'-1\right \|_2.
\end{equation}
\textcolor{orange}{Moreover, the envelope reduction holds:
$\sup_{q\in\mathsf K\cap\mathsf M} M_{U,\varphi}(q)
\ =\
\sup_{q\in\mathsf K\cap\mathsf M} M_{U,\varphi^*}(q),$ 
which recovers the main theorem of \cite{shao2024extreme} in the mean--variance case. }

\end{proposition}

In short, Proposition~\ref{prop:recover-envelope-k2} casts our results as a strict generalization of the convex-envelope approach rather than a competing one: the envelope reduction is recovered without ever convexifying the distortion, as the $U(x)=x$ instance of a single stationarity system that also governs the atomic distortions and the nonconvex, nonsmooth utilities treated in Section~\ref{Examples}.


Mechanistically, the nonconvexity of $\varphi$ is absorbed by the normal-cone correction $z^*\in N_{\mathsf K}(q^*)$ at the quantile level, the solution-level counterpart of convexifying $\varphi$ at the problem level; once $\varphi^*$ is convex this correction is inactive and the system returns the closed-form maximizer and value of Equations \eqref{eq:qstar-envelope}--\eqref{eq:opt-envelope}.

\section{\texorpdfstring{Solution algorithm, closed-form examples, and a real-data application}{Solution algorithm, closed-form examples, and a real-data application}}\label{sec4}

This section turns the optimality theory into a usable solution framework and then puts it to work in a decision. Section~\ref{sec_algorithm} gives the general algorithm, which enforces monotonicity by isotonic projection, and Section~\ref{Examples} uses it to derive closed-form characterizations of the extreme-case value for \changesh{three} representative pairs of utility and distortion, each at a fixed decision. These closed-form inner solves are the building blocks for robust decisions: Section~\ref{sec:genai-ladder} embeds one of them, calibrated on a real serving trace, in a robust capacity-provisioning decision in which the extreme-case value is the inner oracle of an outer min-max problem.

\subsection{Implementation via isotonic projection}\label{sec_algorithm}
Leveraging the optimality conditions established in Section \ref{sec3}, we propose a practical algorithm for solving the general extreme-case DU problem, as detailed in Algorithm \ref{algorithm11}.

\textcolor{orange}{
We illustrate with the minimization problem. By the optimality condition \eqref{optimality condition}, we first construct an unconstrained candidate $R\in\mathcal H=L^2(\mathrm d\nu)$. On the Lebesgue part, it solves $-\lambda_1-2\lambda_2 R(u) \in w(u)\,\partial U(R(u))$ for $\mathrm du$-a.e.\ $u$. For a linear utility, $R_{\lambda_1,\lambda_2}(u):=\left (-\lambda_1-\xi\,w(u)\right )/\left (2\lambda_2\right )$, where $\xi$ denotes the slope of $U$; for a piecewise-linear $U$, after fixing a slope field $\xi(u)$ from the relevant Clarke subdifferential, we set $R(u):= (-\lambda_1-\xi(u)\, w(u) )/ (2\lambda_2 )$ on the Lebesgue part. For the general distortion case in Theorem~\ref{main_thm2_general}, the atomic coordinates of $R$ are understood through the full residual relation $2\lambda_2(R-q^*)=-\rho_\varphi\widehat\xi-(\lambda_1+2\lambda_2q^*)\rho_0$ in $\mathcal H$. }

\textcolor{orange}{
  Next, we project this candidate onto the monotone cone by computing
$$q^*:=\operatorname{Proj}^{\nu}_{\mathsf K}(R) := \arg\min_{q \in \mathsf K} \int_{[0,1]} \left (q(u) - R(u)\right )^2  \mathrm{d}\nu(u).$$ 
{Computationally, the isotonic projection is obtained by iteratively fixing local monotonicity violations:
if two adjacent values are out of order, we combine them into one block and assign the block the $\mathrm d\nu$-weighted average
value. Continuing this procedure eliminates all decreases and yields a globally nondecreasing function.\footnote{ {For details on isotonic projection, see \cite{barlow1972statistical}.}}
}
The associated normal component is
$z^* = -\rho_\varphi\widehat\xi-(\lambda_1+2\lambda_2q^*)\rho_0= 2\lambda_2\left (R - q^*\right ) \in N^{\nu}_{\mathsf K}(q^*)$, 
and can be interpreted as the projection term enforcing monotonicity. Isotonic projection is equivalent to forcing $z$ to satisfy the $L^2(\mathrm d\nu)$ normal-cone complementarity in Theorem~\ref{main_thm2_general}. Figure \ref{Examples for isotonic projection} illustrates three instances of isotonic projection. The red curves are the unconstrained candidates and the blue curves their monotone projections. }

\begin{figure}[htbp]
    \centering
    \begin{subfigure}[t]{0.27\textwidth}
        \centering
        \includegraphics[width=\linewidth]{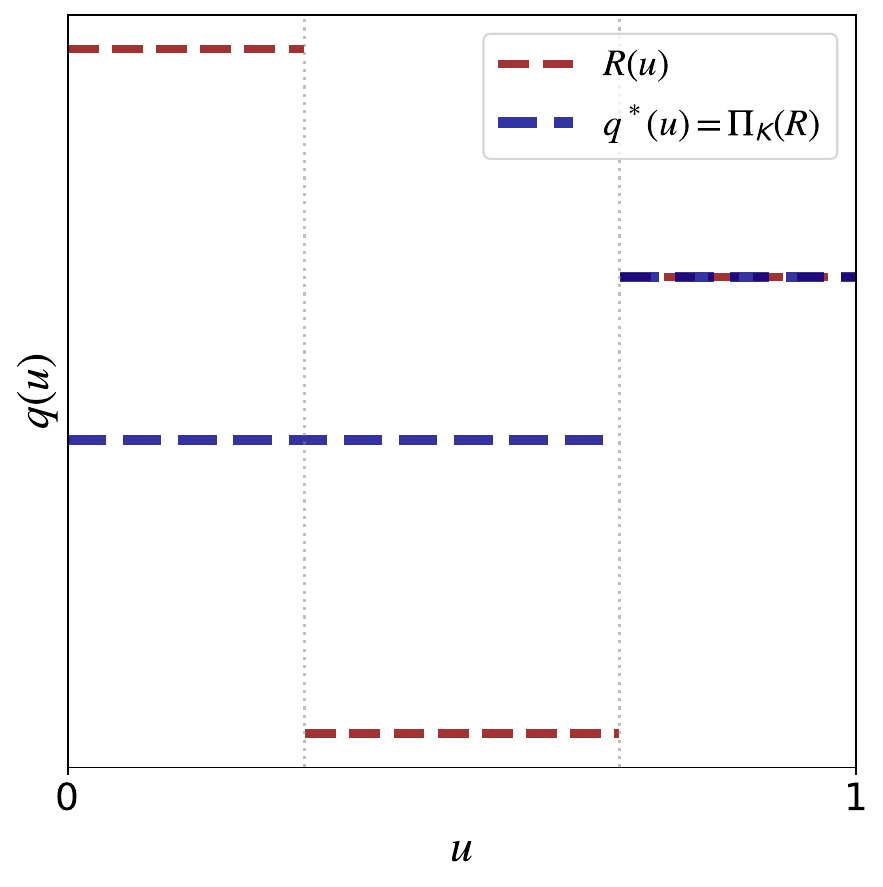}
    \end{subfigure}\hfill
    \begin{subfigure}[t]{0.27\textwidth}
        \centering
        \includegraphics[width=\linewidth]{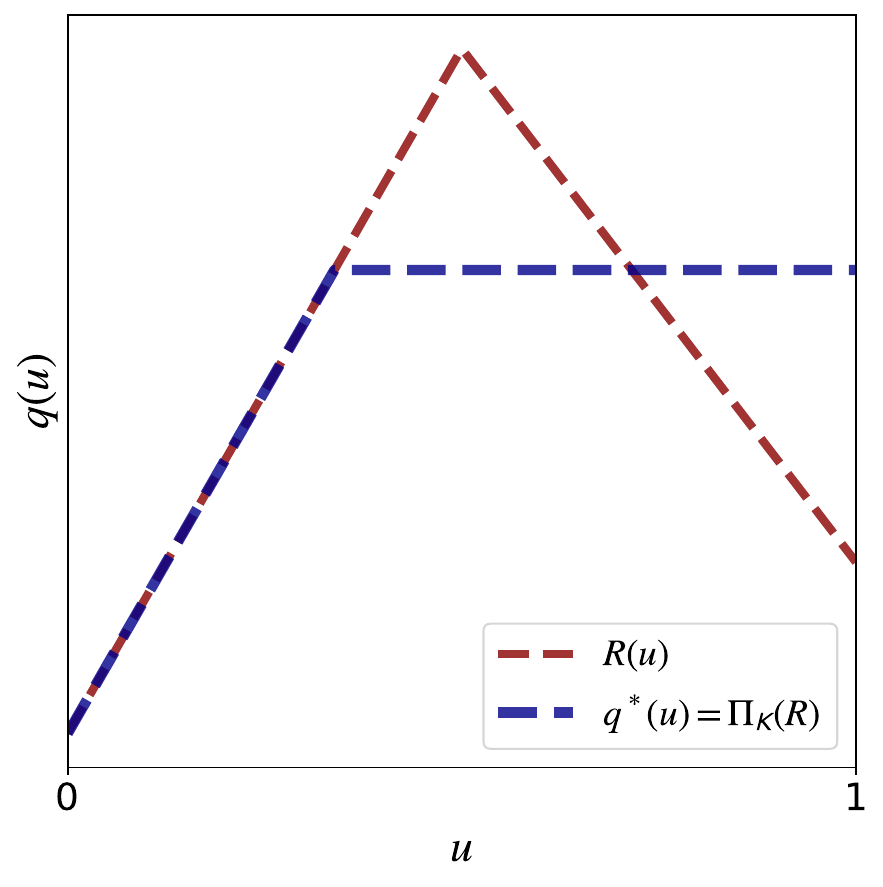}
    \end{subfigure}\hfill
    \begin{subfigure}[t]{0.27\textwidth}
        \centering
        \includegraphics[width=\linewidth]{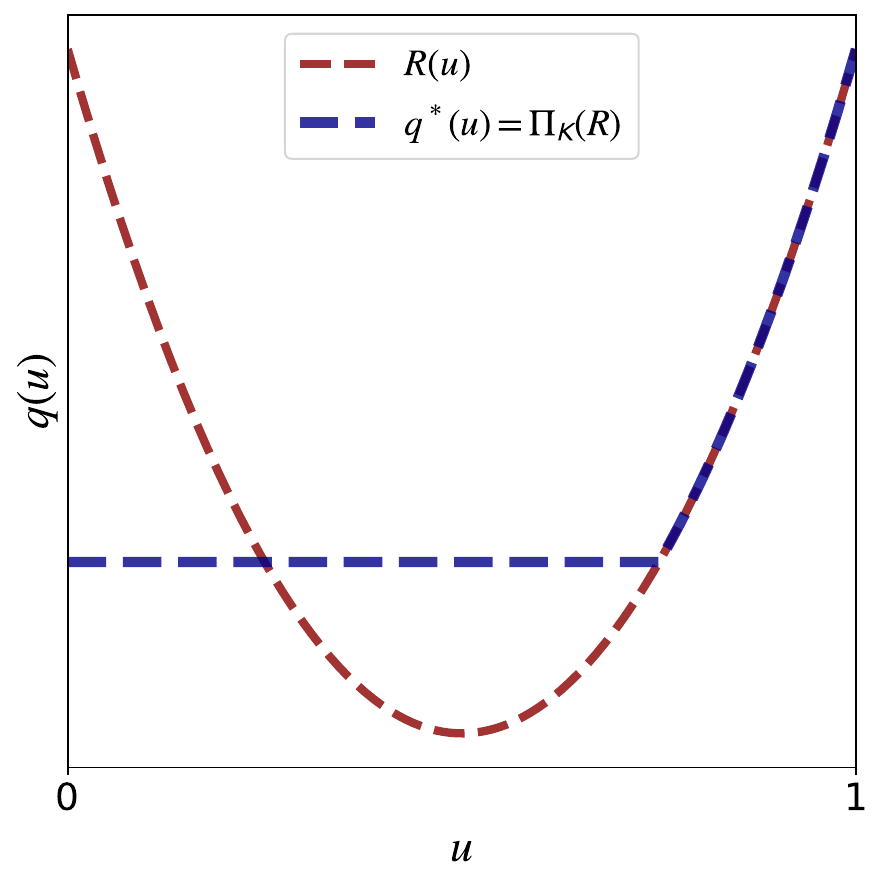}
    \end{subfigure}
    \caption{Three instances of isotonic projection onto the monotone cone: the red curve is the unconstrained candidate $R(u)$ and the blue curve its monotone projection.}
    \label{Examples for isotonic projection}
\end{figure}

{
\renewcommand{\algorithmicrequire}{\textbf{Input:}}
\renewcommand{\algorithmicensure}{\textbf{Output:}}
\begin{algorithm}[htbp!]
\caption{Extreme-case DU values and distributions via isotonic projection.}
\label{algorithm11}
\begin{algorithmic}[1]
\Require Utility $U$, distortion representation $\mathrm d\varphi=w\,\mathrm du+\sum_j m_j\delta_{t_j}$, moments $(\mu,\sigma^2)$, sense $\min/\max$
\State For multipliers $(\lambda_1,\lambda_2)$, construct $R_{\lambda}\in\mathcal H=L^2(\mathrm d\nu)$; on the Lebesgue part it satisfies
\[
0 \in s\, w(u)\,\partial U\!\left (R_{\lambda}(u)\right )+\lambda_1+2\lambda_2 R_{\lambda}(u)
\quad \text{for }\mathrm du\text{-a.e. }u\in(0,1),
\]
where $s=+1$ for minimization and $s=-1$ for maximization, and atomic components are incorporated through the $L^2(\mathrm d\nu)$ residual relation in Theorem~\ref{main_thm2_general}.
\State Enforce monotonicity: $q_{\lambda} \gets \mathrm{Proj}^{\nu}_{\mathsf K}(R_{\lambda})$.
\State Choose $\left (\lambda_1^*,\lambda_2^*\right )$ so that $q_{\lambda^*}$ matches the moments:
$
\int_0^1 q_{\lambda^*}(u)\,\mathrm{d}u=\mu$, $\int_0^1 q^2_{\lambda^*}(u)\,\mathrm{d}u=\mu^2+\sigma^2.
$
\State Output $q^* \gets q_{\lambda^*}$ and recover the optimizer distribution by $
F^*(x):=\sup\!\left \{u\in[0,1]: q^*(u)\le x\right \}$. 
\Ensure Optimal quantile $q^*$ and induced distribution $F^*$
\end{algorithmic}
\end{algorithm}}

Once the monotone quantile $q$ is obtained from $(\lambda_1,\lambda_2)$, we choose the multipliers so that $q$ matches the first two moments. In piecewise-linear settings, isotonic projection makes $q$ piecewise affine in $(\lambda_1,\lambda_2)$, allowing closed-form multipliers; in smooth settings, we solve the two moment equations together with the pointwise stationarity relation defining the unconstrained candidate.

Finally, we recover the optimal CDF by setting $F^*$ to be the generalized inverse of $q^*$. Plateaus of $q^*$ produce atoms of $F^*$; strictly increasing parts produce continuous parts of $F^*$. We verify the complementarity conditions \textcolor{orange}{by computing $z^*=2\lambda_2\left (R-q^*\right )$ }and $\mathcal{J}_{\nu}\!z^*(t)=\int_{[0,t]} z^*(s)\,\mathrm{d}\nu(s)$, and checking that $\mathcal{J}_{\nu}\!z^*(t)\ge0$, $\mathcal{J}_{\nu}\!z^*(1)=0$, and $\int_0^1 \mathcal{J}_{\nu}\!z^*(t)\,\mathrm{d}q^*(t)=0$.

\textcolor{orange}{
The projection step in Algorithm~\ref{algorithm11} supplies the normal-cone component required in Theorem~\ref{main_thm2_general}: if $q^*=\operatorname{Proj}^{\nu}_{\mathsf K}(R)$, then its residual $z^*=2\lambda_2\left (R-q^*\right )$ lies in $N^{\nu}_{\mathsf K}(q^*)$ with $\left \langle z^*,q^*\right \rangle_{\mathcal H}=0$ when $\lambda_2>0$. Therefore, when the constructed $R$, the selectors, and the multipliers also satisfy the stationarity relation and the two moment equations, the resulting $q^*$ satisfies the optimality conditions of Theorem~\ref{main_thm2_general}. Proposition~B.1 in the online appendix states the projection-normal-cone relation formally.
}

\begin{remark}[Computational cost]\label{rem:complexity}
The solution mechanism is computationally inexpensive. On a grid of $n$ quantile levels, the isotonic projection is computed by the pool-adjacent-violators algorithm in $O(n)$ time and memory, and the multipliers are recovered from the moment equalities by solving two equations in the first-two-moment case, or $n$ equations under the first $n$ moments. The unconstrained candidate is available in closed form whenever the utility is piecewise linear, or smooth with an invertible derivative. The overall inner solve therefore scales linearly in $n$ and avoids the semidefinite or conic programs used by relaxation-based approaches, whose size grows with the number of moments and support points. This low and predictable cost makes the extreme-case value and its extremal distribution well suited to serve as an inner oracle that is called repeatedly within a distributionally robust min-max decision loop.
\end{remark}

\subsection{Closed-form examples}\label{Examples}

In this subsection, we apply the proposed algorithm to derive analytical solutions for the optimization problem \eqref{p1}. 
We denote the optimal value-quantile pairs for the maximization and minimization problems by $\left (\Psi_{\max},\,  q_{\max}^*\right )$ and $\left (\Psi_{\min},\,  q_{\min}^*\right )$, respectively.
We present \changesh{three} representative examples characterized by distinct classes of utilities $U$ and distortion functions $\varphi$. For each case, given fixed moment parameters $(\mu, \sigma^2)$, we explicitly construct the optimal quantile function $q^*$ that satisfies the variational inequality system \eqref{eq:VI}, thus providing exact characterizations for both the minimization and maximization problems. Throughout, a closed-form solution refers to an explicit finite-dimensional characterization rather than a solution of the infinite-dimensional functional program or a large-scale conic relaxation: the extremal distribution is a finite staircase or an explicit density, and the optimal value together with the multipliers is pinned down by an explicit system that requires at most an isotonic projection and a scalar root-finding step.

\changesh{The framework imposes no structural restriction on the pair $(U,\varphi)$: the utility may be smooth or nonsmooth, and the distortion may be absolutely continuous or carry atoms. In the main text we concentrate on the cases that are most salient for operations and least accessible to existing theory, those in which the extremal distribution is a low-dimensional discrete law, the regime produced by a piecewise-linear utility or by a distortion with atoms. The complementary smooth case, in which a smooth utility and an absolutely continuous distortion yield a continuous worst-case distribution, is handled by the same machinery; we work out one such instance, a Dual--Power spectral distortion with quadratic utility, in Appendix~B.3, where the worst case also exhibits a finite but unattained supremum.}

\changesh{The three main-text examples span this range. Example~1 extends the classical Scarf bound from an average criterion to an RVaR tail window; Example~2 treats a GlueVaR distortion that carries atoms, paired with a reward--penalty utility; and Example~3 crosses into the genuinely nonconvex and nonconcave regime under a CVaR distortion, the case that no prior moment method covers and that motivates the framework. Each example is a closed-form characterization of the extreme-case value for a fixed decision, and Section~\ref{sec:genai-ladder} illustrates how these inner characterizations support a real robust decision application.}

\subsection*{Example 1: Scarf bound with RVaR distortion}
The shortfall payoff $U(x)=\left (x-a\right )_+$ is a canonical nonsmooth function in moment-based DRO, isolating upper-tail outcomes above a threshold $a$. In the classical Scarf bound, the objective is the extreme-case expectation and the worst-case distribution is two-point \citep{scarf1958min,jagannathan1977minimax,kuhn2025distributionally,jiang2025dual}. We extend this bound from identity distortion to RVaR distortion and obtain explicit closed forms. RVaR upgrades Scarf-type shortfall control from an average criterion to a tunable focus on the tail window indexed by $(\alpha,\beta)$ \citep{li2018worst}.

\begin{proposition}[Scarf bound under RVaR]\label{prop:Scarf-RVaR}
Define  $U(x):=\left (x-a\right )_+$ for a threshold  $a\in\mathbb R$ and  $w_{\alpha,\beta}(u):=\mathbf 1_{(\alpha,\beta]}(u)/(\beta-\alpha)$ for  two confidence levels $0\le \alpha<\beta\le 1$.  For any $\nu\in(0,1)$, let  
$$L_\nu:=\mu-\sigma\sqrt{\frac{1-\nu}{\nu}},\
U_\nu:=\mu+\sigma\sqrt{\frac{\nu}{1-\nu}}. $$
Then the following conclusions hold.

\begin{enumerate}[(a)]
\item 
$\Psi_{\min}=\max\!\left \{L_\beta-a, 0\right \}$,  which is attained by the two-level quantile
$q^*_{\min}(u)
=
L_\beta\,\mathbf 1_{\left \{u\le\beta\right \}}
+
U_\beta\mathbf 1_{\left \{u>\beta\right \}}.$

\item Define $T_{\alpha,\beta}:=\mu+\sigma \frac{-2\alpha^2+3\alpha-\beta}{2\left (1-\alpha\right )\sqrt{\alpha\left (1-\alpha\right )}}$, which obviously satisfies  $T_{\alpha,\beta}\le U_\alpha<U_\beta$. For $T_{\alpha,\beta}<a<U_\beta$, let $t^*_{\alpha,\beta}$ be the unique positive root of
\begin{equation}\label{eq:quartic}
(1-\beta)\,t^4+(3-2\beta)\,t^2+\frac{2\left (\mu-a\right )}{\sigma}\,t-\beta=0,
\end{equation}
and set $p^*:={(t^*_{\alpha,\beta})^2}/(1+(t^*_{\alpha,\beta})^2)$.
Then $\Psi_{\max}$ admits the following
 closed-form solution:
\[
\Psi_{\max}=
\begin{cases}
U_\alpha-a,&  a\le T_{\alpha,\beta},\\[2pt]
\frac{\beta-p^*}{\beta-\alpha}\left (\mu-a+\sigma t_{\alpha,\beta}^*\right ),&
 T_{\alpha,\beta}<a<U_\beta,\\[2pt]
0,& a\ge U_\beta, 
\end{cases}
\]
Moreover, the associated maximizers take the following forms:
\begin{displaymath}
q^*_{\max}(u)=
\begin{cases}
L_\alpha \mathbf 1_{\left \{u\le\alpha\right \}}+U_\alpha \mathbf 1_{\left \{u>\alpha\right \}}, & a\le T_{\alpha,\beta},\\
L_{p^*} \mathbf{1}_{\left \{0\le u\le p^*\right \}}+U_{p^*}\mathbf 1_{\left\{ p^*<u\le 1\right \}}, &  T_{\alpha,\beta}<a<U_\beta,\\
q^*_{\min}(u), & a\ge U_\beta,    
\end{cases}\quad\text{for } u\in[0,1]. 
\end{displaymath}

\item When $\alpha=0$ and $\beta=1$, $w_{\alpha,\beta}(u)=1$ for all $u\in[0,1]$, which exactly corresponds to the expectation. In this case,  we obtain 
$\displaystyle \lim_{\alpha \downarrow 0,\,  \beta \uparrow 1} T_{\alpha,\beta}=-\infty$ and $\displaystyle \lim_{\beta\uparrow 1} U_\beta=+\infty$, and Equation \eqref{eq:quartic} reduces to a quadratic function with 
$$t^*_{0,1}=\frac{\sqrt{\left (\mu-a\right )^2+\sigma^2}-\left (\mu-a\right )}{\sigma},$$ which further yields 
$${\sup_{q\in\mathsf K\cap\mathsf M}\Psi(q)}=\frac12\left (\sqrt{\left (a-\mu\right )^2+\sigma^2}-\left (a-\mu\right )\right ).$$
Moreover, the maximizer takes the two-point form with values $a\pm\sqrt{\left (\mu-a\right )^2+\sigma^2}$,  recovering the classical Scarf bound.
\end{enumerate}
\end{proposition}

Figure \ref{fig:scarf_quantiles} shows the best-case quantile and the two nontrivial worst-case quantiles, together with the corresponding unconstrained candidates $R$ before projection. Panel (a) gives the best case, a two-point distribution jumping at $\beta$ as in Proposition~\ref{prop:Scarf-RVaR}(a), while panels (b) and (c) give the two nontrivial worst-case regimes. Figure \ref{fig: Psi via a} plots the worst-case risk $\Psi_{\max}$ against the threshold $a$ as a nonincreasing continuous function with a threshold-driven regime structure determined by $a$ and by $\alpha,\beta$ in the RVaR weight.

\begin{figure}[htbp]
    \centering
    \includegraphics[width=0.88\textwidth,height=0.28\textheight]{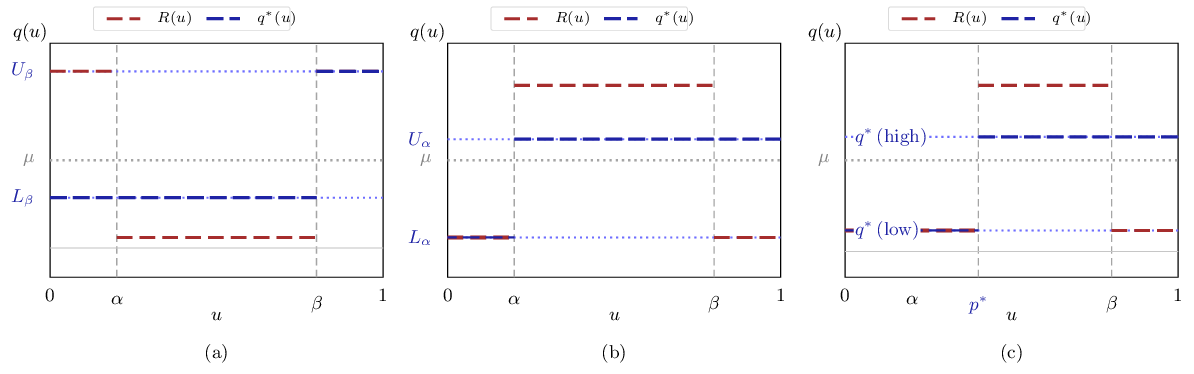}
    \caption{Extreme quantiles for the RVaR-distorted Scarf example of Proposition~\ref{prop:Scarf-RVaR}. Panel (a) gives the best-case quantile $q^*_{\min}(u)=L_\beta\mathbf 1_{\{u\le\beta\}}+U_\beta\mathbf 1_{\{u>\beta\}}$ and unconstrained candidate $R(u)$. Panels (b) and (c) give the low-threshold and interior worst-case quantiles, jumping at $\alpha$ and $p^*$.}
    \label{fig:scarf_quantiles}
\end{figure}
For low thresholds $a\le T_{\alpha,\beta}$, the worst-case RVaR is attained by a two-point distribution with breakpoint $\alpha$; for sufficiently large $a\ge U_{\beta}$, tail losses vanish and $\Psi_{\max}=0$. In the intermediate regime, the binding distortion induces an interior two-point extremal distribution whose mass split is determined by a one-dimensional equation, so $\Psi_{\max}$ depends nonlinearly on $a$ rather than linearly as under expectation-type objectives.

Proposition~\ref{prop:Scarf-RVaR}(c) recovers the classical degenerate case \citep{scarf1958min,jagannathan1977minimax}. When $\beta\downarrow\alpha$, RVaR degenerates to VaR, and when $\beta\uparrow 1$, it degenerates to CVaR, so the analysis recovers the Scarf bound in the VaR and CVaR cases \citep{cai2024worst}. 
Consequently, the distortion modifies the worst case only through the locations and weights
of the extremal atoms, preserving the tractability and interpretability of Scarf-type bounds
while extending them to a broad class of  tail-sensitive risk measures.
\begin{figure}[htbp]
    \centering
    \includegraphics[width=0.4\textwidth]{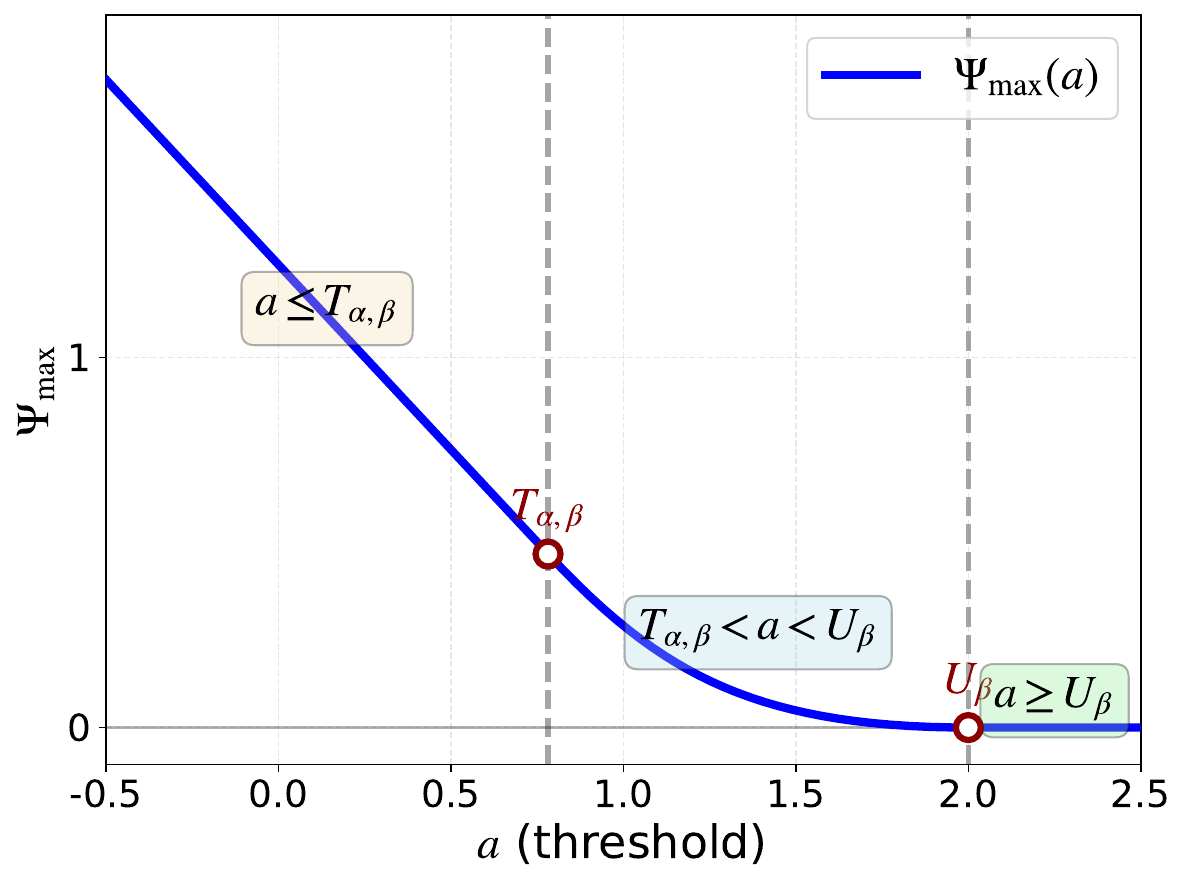}
    \caption{Worst-case RVaR-distorted value $\Psi_{\max}(a)$ as a function of the threshold $a$ in the Scarf example of Proposition~\ref{prop:Scarf-RVaR}. The function is continuous and nonincreasing: affine with slope $-1$ for $a\le T_{\alpha,\beta}$, nonlinear for $T_{\alpha,\beta}<a<U_\beta$, and zero for $a\ge U_\beta$.}
    \label{fig: Psi via a}
\end{figure}

\subsection*{Example 2: Reward--penalty utility with GlueVaR}

We next apply Theorem \ref{main_thm2_general} to general DRMs that need not be absolutely continuous and may be nonconvex, extending the conclusions of \cite{cornilly2018upper} and \cite{li2018closed}.
A prominent example is GlueVaR, which combines an absolutely continuous component
with point masses at selected quantile levels. The distortion function of GlueVaR takes the following form: 
\begin{equation}\label{glue_var}
\varphi(u)=\gamma\,\mathbf 1_{[\alpha,\beta]}(u)+\frac{\left (1-\gamma\right )u+\gamma-\beta}{1-\beta}\,\mathbf 1_{(\beta,1]}(u),\quad 0<\alpha<\beta<1,\ \gamma\in(0,1)
\end{equation}
 It has been widely studied in risk management as a flexible framework that interpolates
between VaR- and CVaR-type measures and allows practitioners to assign explicit importance
to specific tail events \citep{belles2014beyond,belles2016attitudes}. To connect this risk measure with a practically interpretable payoff, we consider the piecewise-linear utility $U(x)=x+\kappa\left (x-a\right )_+$, which prices the baseline outcome linearly but imposes an additional penalty once the loss exceeds a target threshold
$a$, as in a service-level penalty contract. Increasing $\kappa$ explicitly amplifies the contribution of exceedance events beyond $a$, aligning naturally with GlueVaR’s ability to place extra emphasis on selected tail quantiles via point masses \citep{cai2025conditional}.

\begin{proposition}[Reward-penalty utility with GlueVaR]\label{prop:gluevar-extrema}
Define $U(x)=x+\kappa\left (x-a\right )_+$ with  $\kappa>0$ and let  the distortion function be specified in Equation \eqref{glue_var}. Denote {the constants}
\begin{displaymath}
D:=\sqrt{\alpha+\frac{\left (\gamma-\beta+\alpha\right )^2}{\beta-\alpha}+\frac{\left (\beta-\gamma\right )^2}{1-\beta}},\
D_\kappa:=\frac{\sqrt{\alpha B_\kappa^2+\left (\beta-\alpha\right )\left (B_\kappa-\underline{\rho}\right )^2+\left (1-\beta\right )\left (B_\kappa-1\right )^2}}{B_\kappa},
\end{displaymath}
where 
$$\underline{\rho}:=\frac{\gamma\left (1-\beta\right )}{\left (1+\kappa\right )\left (1-\gamma\right )\left (\beta-\alpha\right )},\ B_\kappa:=1-\beta+\left  (\beta-\alpha\right)\underline{\rho}.$$
Then  the optimal quantile functions and objective values are determined by the following cases. 

\noindent
	 {
(a1) {For the maximization problem,} if $0<\gamma<\frac{\beta-\alpha}{1-\alpha}$, then $q^*_{\max}$ is piecewise constant and admits a closed-form solution.
For $s\in[\beta,1]$, define
\begin{align*}
C(s)&:=\gamma+\frac{1-\gamma}{1-\beta}\left (s-\beta+\left (1+\kappa\right )\left (1-s\right )\right ),\\
\Delta(s)&:=\frac{\gamma^2}{\beta-\alpha}+\left (\frac{1-\gamma}{1-\beta}\right )^2\left (s-\beta+\left (1+\kappa\right )^2\left (1-s\right )\right )-C^2(s),
\end{align*}
and set the two target thresholds
\begin{align}
a_{\min}&:=\mu-\frac{\sigma}{2\sqrt{\Delta(\beta)}}\left (2\left (1+\kappa\left (1-\gamma\right )\right )-\frac{1-\gamma}{1-\beta}\left (2+\kappa\right )\right ),\label{amin}\\
a_{\max}&:=\mu-\frac{\sigma}{2\sqrt{\Delta(1)}}\left (2-\frac{1-\gamma}{1-\beta}\left (2+\kappa\right )\right ).\label{amax}
\end{align}
Except when $a\in[a_{\min},a_{\max}]$, the optimizer has three levels with fixed breakpoints at $\alpha$ and $\beta$.
In the intermediate regime $a\in[a_{\min},a_{\max}]$, the optimizer has four levels: besides $\alpha$ and $\beta$,
an additional breakpoint $s^*\in(\beta,1)$ arises endogenously.}

\smallskip

\begin{enumerate}[(i)]{
\item If $a\le \mu-\frac{\sigma}{D}$, then the maximizer is a three-level quantile:
\begin{equation*}
   q^*_{\max}(u)= \underbrace{\left(\mu-\frac{\sigma}{D}\right )}_{c_1} \mathbf 1_{\left \{u<\alpha\right \}}+ \underbrace{\left (\mu+\frac{\sigma\left (\gamma-\beta+\alpha\right )}{D\left (\beta-\alpha\right )}\right )}_{c_2} \mathbf 1_{\{\alpha\le u<\beta\}}+ \underbrace{\left (\mu+\frac{\sigma \left (\beta-\gamma\right )}{D\left (1-\beta\right )}\right )}_{c_3}\mathbf 1_{\{\beta\le u\le 1\}},\ u\in[0,1], 
\end{equation*}
and 
the associated optimal value is $\gamma  U(c_2)+\left (1-\gamma\right )U(c_3)$. 

\item If $\mu-\frac{\sigma}{D}<a<\mu-\frac{\sigma}{D_\kappa}$, the maximizer is a three-level quantile:
\begin{equation*}
   q^*_{\max}(u)=\underbrace{a}_{c_1} \mathbf 1_{\{u<\alpha\}}+ \underbrace{\left (a+ \frac{(\mu-a)\rho(a)}{1-\beta+\left (\beta-\alpha\right )\rho(a)}\right )}_{c_2}\mathbf 1_{\left \{\alpha\le u<\beta\right \}}+ \underbrace{\left (a+\frac{\mu-a}{1-\beta+\left (\beta-\alpha\right )\rho(a)}\right )}_{c_3}\mathbf 1_{\left \{\beta\le u\le 1\right \}}, 
\end{equation*} for $u\in[0,1]$, 
  where $\rho(a)>0$ solves the equation of $\rho$: $$\left (\left (\beta-\alpha\right )\rho^2+(1-\beta)\right )\left (\mu-a\right )^2=\left (\left (\mu-a\right )^2+\sigma^2\right )\left ((1-\beta)+\left (\beta-\alpha\right )\rho\right )^2$$ and additionally satisfies $\underline{\rho}\le \rho(a)\le \left (1+\kappa\right )\underline{\rho}$. Moreover, the associated  optimal value is $\gamma  U(c_2)+\left (1-\gamma\right )U(c_3)$. 

\item If $\mu-\frac{\sigma}{D_\kappa}\le a<a_{\min}$, then the maximizer is a three-level quantile:
\begin{equation*}
   q^*_{\max}(u)=\underbrace{\left (\mu-\frac{\sigma}{D_\kappa}\right )}_{c_1}\mathbf 1_{\left \{u<\alpha\right \}}+\underbrace{\left (\mu+\left (\frac{\underline{\rho}}{B_\kappa}-1\right )\frac{\sigma}{D_\kappa}\right )}_{c_2}\,\mathbf 1_{\left \{\alpha\le u<\beta\right \}}+ \underbrace{\left (\mu+\left (\frac{1}{B_\kappa}-1\right )\frac{\sigma}{D_\kappa}\right )}_{c_3}\,\mathbf 1_{\left \{\beta\le u\le 1\right \}},
\end{equation*}for $u\in[0,1]$, 
 and the  associated optimal value is 
$\Psi_{\max}=\gamma  U(c_2)+\left (1-\gamma\right )U(c_3)$.

\item If $a_{\min}\le a\le a_{\max}$, define
\begin{align*}
c_1(s)&:=\mu-C(s)\frac{\sigma}{\sqrt{\Delta(s)}},\quad c_2(s):=c_1(s)+\frac{\gamma}{\beta-\alpha}\frac{\sigma}{\sqrt{\Delta(s)}},\\ c_3(s)&:=c_1(s)+\frac{1-\gamma}{1-\beta}\frac{\sigma}{\sqrt{\Delta(s)}},\quad c_4(s):=c_1(s)+\frac{(1-\gamma)(1+\kappa)}{1-\beta}\frac{\sigma}{\sqrt{\Delta(s)}}
\end{align*}
Then the maximizer is a four-level quantile
\begin{equation*}
    q^*_{\max}(u)=c_1(s^*)\mathbf 1_{\{u<\alpha\}}+c_2(s^*)\mathbf 1_{\{\alpha\le u<\beta\}}+c_3(s^*)\mathbf 1_{\{\beta\le u< s^*\}}+c_4(s^*)\mathbf 1_{\{s^*\le u\le 1\}},
\end{equation*}
where $s^*\in(\beta,1)$ satisfies
\begin{equation}\label{eq:s}
    \mu- a =\frac{\sigma}{2\sqrt{\Delta(s^*)}}\big (2C(s^*)-\frac{1-\gamma}{1-\beta}\left (2+\kappa\right )\big ).
\end{equation} 
The  associated optimal value is 
$$\Psi_{\max}=\gamma U(c_2(s^*))+\left (1-\gamma\right )\left (\frac{s^*-\beta}{1-\beta}U(c_3(s^*))+\frac{1-s^*}{1-\beta}U(c_4(s^*)) \right ).$$ }

\item If $a\ge a_{\max}$, then the maximizer returns to the three-level form in \textnormal{(i)} with $(c_1,c_2,c_3)$ given by $$c_1:=\mu-\frac{\sigma}{D},\quad c_2:=\mu+\frac{\gamma-\beta+\alpha}{\beta-\alpha}\frac{\sigma}{D},\quad c_3:=\mu+\frac{\beta-\gamma}{1-\beta}\frac{\sigma}{D},$$ and the optimal value is $\Psi_{\max}=\gamma\, U(c_2)+\left (1-\gamma\right )U(c_3)$.

\end{enumerate}

\noindent
{(a2)  For the maximization problem,} if $\gamma\ge\frac{\beta-\alpha}{1-\alpha}$, denote
$L_\alpha:=\mu-\sigma\sqrt{\frac{1-\alpha}{\alpha}}$,
$L_\beta:=\mu-\sigma\sqrt{\frac{1-\beta}{\beta}}$, 
$U_\alpha:=\mu+\sigma\sqrt{\frac{\alpha}{1-\alpha}}$, 
$U_\beta:=\mu+\sigma\sqrt{\frac{\beta}{1-\beta}}$. then the isotonic projection pools the last two blocks and the maximizer degenerates to a two-level quantile $q^*_{\max}(u)=L_\alpha\,\mathbf 1_{\{0\le u<\alpha\}}+U_\alpha\,\mathbf 1_{\{\alpha\le u\le 1\}}$. Consequently, the optimal value depends only on the position of $a$ relative to $U_\alpha$: if $a\le U_\alpha$, then $\Psi_{\max}=U(U_\alpha)=\left (1+\kappa\right )U_\alpha-\kappa a$, and if $a\ge U_\alpha$, then $\Psi_{\max}=U(U_\alpha)=U_\alpha$.

\noindent 
{(b) For the minimization problem,} if $\gamma\ge \frac{U(U_{\alpha})-U(\mu)}{U(U_{\alpha})-U(L_\alpha)}$, then the infimum is attained with minimizer 
$q_{\min}^*(u)=L_\alpha\,\mathbf 1_{\{u\le \alpha\}}+U_\alpha\,\mathbf 1_{\{u>\alpha\}},$
and the minimum value equals
$\Psi_{\min}=\gamma\,U(L_\alpha)+(1-\gamma)\,U(U_\alpha).$ If $\gamma < \frac{U(U_{\alpha})-U(\mu)}{U(U_{\alpha})-U(L_\alpha)}$, then the infimum $U(\mu)$ is not attained.

\end{proposition}

Proposition~\ref{prop:gluevar-extrema} characterizes the extreme values of the GlueVaR-distorted objective under the reward--penalty utility. For the maximization part, the atom at $\alpha$ and the tail weight beyond $\beta$ produce a richer extremal geometry, with the maximizer determined by the kink location $a$ through the subdifferential of $U$. When $0<\gamma<(\beta-\alpha)/(1-\alpha)$, it is a three-level quantile with breakpoints at $\alpha,\beta$ in regimes (i)--(iii), becomes genuinely four-level on the intermediate range $a_{\min}\le a\le a_{\max}$ where the tail block splits at an endogenous breakpoint $s^*\in(\beta,1)$ solving Equation \eqref{eq:s}, and flattens once $a\ge a_{\max}$ as the penalty turns inactive; when $(\beta-\alpha)/(1-\alpha)\le\gamma\le1$ it degenerates to two levels as monotonicity pools the last two blocks.
For the minimization part, the infimum exhibits a Scarf-type two-point extremal structure.

Figure~\ref{gluevar} plots the worst-case value $\Psi_{\max}$ as a nonincreasing continuous function of $a$ when $0<\gamma<\frac{\beta-\alpha}{1-\alpha}$. It illustrates the five cases in Proposition \ref{prop:gluevar-extrema}(a), highlighting the regime switches at $\mu-\sigma/D$, $\mu-\sigma/D_\kappa$, $a_{\min}$, and $a_{\max}$, and the plateau beyond $a_{\max}$. This example shows the value of the projection framework: it handles a DRM with atoms and a nonsmooth, nontrivial utility, yet still yields explicit extremal quantiles and tractable values in terms of $(\alpha,\beta,\gamma,\kappa,\mu,\sigma)$. Proposition~\ref{prop:gluevar-extrema} reduces to the worst-case GlueVaR under the identity utility by setting $\kappa=0$ \citep{shao2024extreme}, and recovers the CVaR special case of the same reward--penalty utility by taking $\alpha=\beta$ and $\gamma=0$ \citep{cai2025conditional}.\footnote{{When $\alpha=\beta$ and $\gamma=0$, the standardized thresholds $(a-\mu)/\sigma$ in Equations \eqref{amin} and \eqref{amax} reduce to $a_{\min}=\frac{1+(1+\kappa)(2\alpha-1)}{2(1+\kappa)\sqrt{\alpha(1-\alpha)}}$ and $a_{\max}=\frac{\kappa+2\alpha}{2\sqrt{\alpha(1-\alpha)}}$ in \cite{cai2025conditional}. Cases in (i),(ii),(iii),(v) reduce to two-level quantiles, while the case (iv) reduces to a three-level quantile due to the combination of $\alpha$ and $\beta$. }}

\begin{figure}[htbp]
    \centering
    \includegraphics[width=0.5\textwidth]{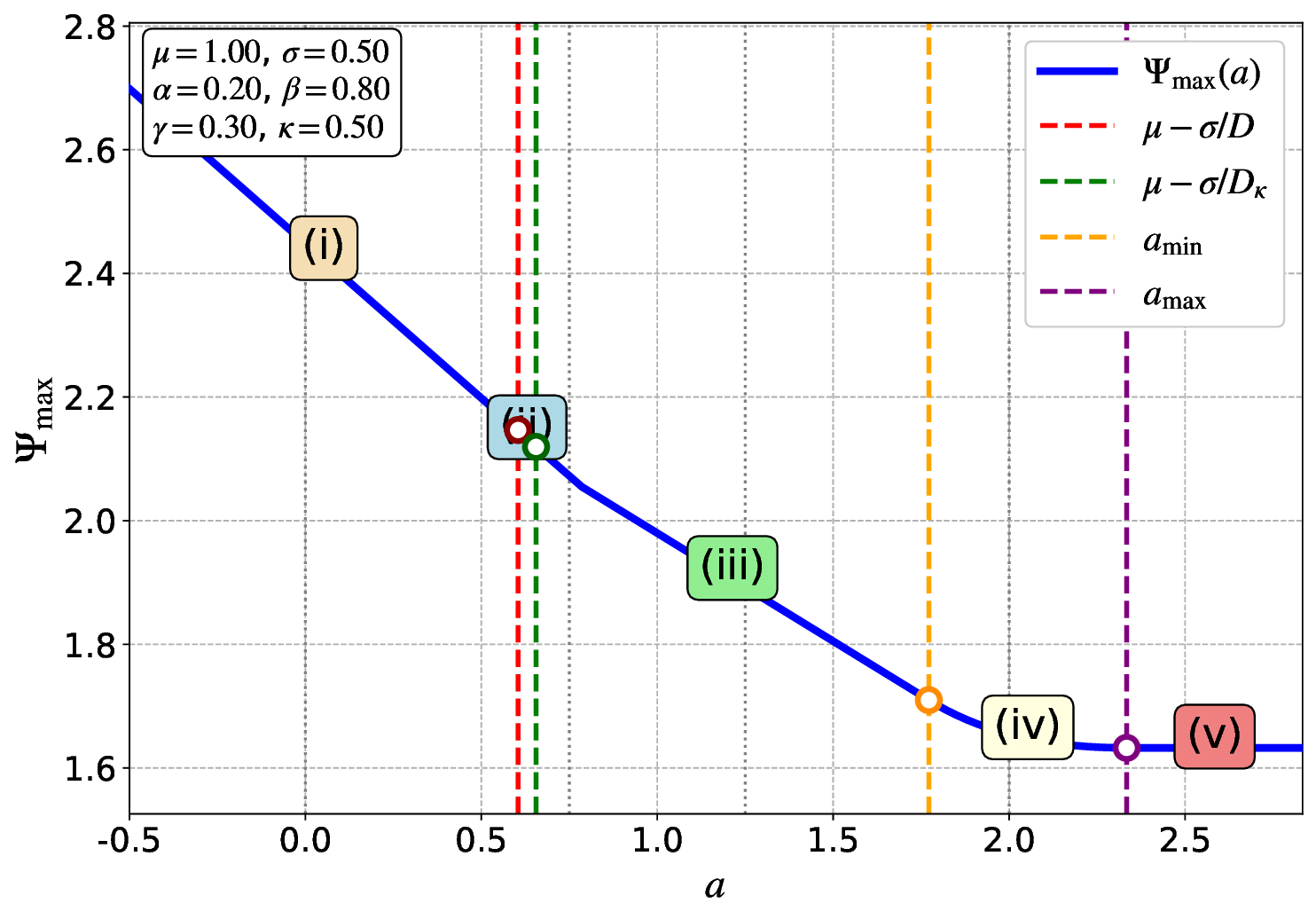}
    \caption{Worst-case GlueVaR-distorted value $\Psi_{\max}(a)$ as a function of the threshold $a$ for the reward--penalty example of Proposition~\ref{prop:gluevar-extrema}, with $0<\gamma<(\beta-\alpha)/(1-\alpha)$. The curve is continuous, nonincreasing, and follows the five cases of \changesh{part~(a)}, with switches at $\mu-\sigma/D$, $\mu-\sigma/D_\kappa$, $a_{\min}$, $a_{\max}$, and a \changesh{final} plateau.}
    \label{gluevar}
\end{figure}

\subsection*{Example 3: Robust design of a capped incentive contract}
\changesh{Examples~1 and~2} evaluate a convex utility under a nontrivial distortion. We now study a contract-design problem that is genuinely neither convex nor concave in its payoff and is evaluated through a nontrivial distortion at the same time, so that the worst case is an RDU functional of the outcome distribution rather than a linear one. The classical moment duality that solves a plain expected payoff does not apply, and the quantile projection of Section~\ref{sec3} is what delivers the closed form. A principal rewards an agent for a performance outcome $X$ through a capped bonus
\begin{equation}\label{eq:capped-U}
U(x)=\min\!\left (\kappa\left (x-a\right )_+, B\right ),\quad \kappa>0,\ B>0,
\end{equation}
which pays at a marginal rate $\kappa$ once the outcome clears a target $a$ and caps the total bonus at a budget $B$. With $\bar a:=a+B/\kappa$ the outcome at which the cap is first reached, the bonus is the vertical-spread payoff $U(x)=\kappa\left (x-a\right )_+-\kappa\left (x-\bar a\right )_+$, which is nondecreasing and Lipschitz with a convex kink at $a$ and a concave kink at $\bar a$, hence neither convex nor concave. This capped, limited-loss payoff is the worst-case object studied in robust reinsurance under moment information \citep{cai2024worst}, which our distortion-aware analysis extends. The principal does not budget the bonus at its mean but provisions it by CVaR of the bonus at level $\alpha$, the average bonus over the most expensive $1-\alpha$ fraction of outcome scenarios, which for an increasing payoff are the high-output ones. The CVaR is the tail case $\beta\uparrow 1$ of the RVaR distortion. Only the first two moments $(\mu,\sigma^2)$ of the outcome are known, and the limit $\alpha\downarrow 0$ removes the distortion and returns the plain expected bonus.

\begin{proposition}[Capped incentive contract under CVaR provisioning]\label{prop:capped}
Let $U$ be given by Equation \eqref{eq:capped-U} with $\bar a=a+B/\kappa$, let $\varphi(u)=\max\!\left \{0,\left (u-\alpha\right )/\left (1-\alpha\right )\right \}$ for some $\alpha\in(0,1)$, and write
\[
U_\alpha=\mu+\sigma\sqrt{\frac{\alpha}{1-\alpha}},\quad
T_\alpha=\mu+\sigma\,\frac{2\alpha-1}{2\sqrt{\alpha\left (1-\alpha\right )}}.
\]
Let $S(a)$ be the uncapped worst case of Example~1 for the ramp $\kappa(x-a)_+$ under the same distortion,
\begin{displaymath}
S(a)=\begin{cases}\kappa\left (U_\alpha-a\right ),& a\le T_\alpha,\\[3pt]\dfrac{\kappa}{2\left (1-\alpha\right )}\left (\sqrt{\left (a-\mu\right )^2+\sigma^2}-\left (a-\mu\right )\right ),& a> T_\alpha.\end{cases}
\end{displaymath}
Then the worst-case provisioned bonus is
\begin{equation}\label{eq:capped-sup}
\Psi_{\max}=\begin{cases}
B,& \bar a\le U_\alpha,\\[2pt]
\frac{B \sigma^2}{\left (1-\alpha\right )\left (\sigma^2+\left (\bar a-\mu\right )^2\right )},& \bar a> U_\alpha\ \text{ and }\ \sqrt{\left (a-\mu\right )^2+\sigma^2}\ge B/\kappa,\\[4pt]
S(a),& \bar a> U_\alpha\ \text{ and }\ \sqrt{\left (a-\mu\right )^2+\sigma^2}< B/\kappa,
\end{cases}
\end{equation}
attained by a two-point distribution. The best-case provisioned bonus satisfies $\Psi_{\min}=0$ for $a\ge\mu$, attained by a two-point distribution supported on $(-\infty,a]$ when $a>\mu$ and approached when $a=\mu$, and $\Psi_{\min}>0$ for $a<\mu$.
\end{proposition}

Three readings make the result operational. First, distortion makes the problem nonclassical: because $M_{U,\varphi}$ is an RDU functional, the worst case cannot be obtained from the quadratic-majorant duality for the plain expected bonus, so the quantile projection of Section~\ref{sec3} is required. As $\alpha\downarrow 0$, $U_\alpha$ collapses to $\mu$ and $S(a)$ to ${\kappa}\left (\sqrt{\left (a-\mu\right )^2+\sigma^2}-\left (a-\mu\right )\right )/2$, recovering the plain expected-bonus bound; hence the undistorted classical problem is nested in our characterization. Second, the cap remains a hard budget guarantee, since $\Psi_{\max}\le B$ for every distribution, and CVaR provisioning widens the budget-saturation regime from $\bar a\le\mu$ to $\bar a\le U_\alpha$. Third, when the cap is slack, the worst case coincides with the uncapped value $S(a)$ of Example~1; the only new behavior is the binding regime, in which the worst-case reserve $B\sigma^2/\left (\left (1-\alpha\right )\left (\sigma^2+\left (\bar a-\mu\right )^2\right )\right )$ rises with $\alpha$ and $\sigma$ and falls as the target or cap reach $\bar a$ is pushed out.

The examples above characterize the extreme-case value for a fixed decision. We now use this characterization as the inner problem of an outer min-max decision, in which a decision maker optimizes a variable against the worst case over the moment ambiguity set.

\subsection{\texorpdfstring{Application: robust capacity provisioning for generative artificial intelligence inference}{Application: robust capacity provisioning for generative artificial intelligence inference}}\label{sec:genai-ladder}

\textbf{The decision.} An inference platform serves a large-language-model API to enterprise customers and must commit accelerator capacity for the coming day before that day's demand is known, through reserved-capacity instruments such as fixed-window GPU reservations or pre-purchased provisioned throughput. The governing service-level agreement (SLA) reports realized load at a few named tail percentiles---in our data the 90th, 95th, and 99th---and pays the customer a service credit whenever a reported percentile overruns the reserved capacity. The credit schedule is \emph{tiered and capped}: it pays nothing until an overrun threshold, then steps up band by band as the breach deepens, and saturates at a contractual ceiling. Reserving too much strands expensive accelerators; reserving too little triggers credits; and because token traffic is heavily bursty, the credit bill is dominated by deep tail events. The platform's daily problem, repeated across a fleet of models and regions, is to reserve the \emph{least} capacity whose worst-case credit liability stays within a fixed budget---the error budget of site-reliability practice. These contractual features are exactly the primitives of our model: the tiered, capped credit is the multi-band utility~\eqref{eq:gl-U}; the named-percentile reporting is the atomic distortion; and the liability budget is the constraint $\Psi_{\max}(c)\le B_0$ in~\eqref{eq:gl-outer}.

\textbf{Why standard sizing tools fall short.} Planners today size capacity from a single distribution: a Gaussian or lognormal fit to recent load, training-window empirical percentiles, or a per-percentile Cantelli (mean--variance) safety margin. Each commits to one law inside the set of distributions consistent with the data, so each can be undercut on precisely the tail it did not anticipate. Worse, the decision carries three features that defeat classical worst-case tools at once. The credit is tiered and capped, so liability is neither convex nor concave and the convex-envelope reductions that handle linear or concave payoffs do not apply \citep{shao2024extreme, pesenti2024optimizing, cai2025distributionally}. It is charged on a few named percentiles, so the distortion is atomic, outside the absolutely continuous theory. And the load is known only through its mean and variance---least reliably in the deep tail that drives the bill. The robust response is to size capacity against the worst distribution matching the moments: an extreme-case distorted utility with a nonconvex payoff and an atomic distortion. This is exactly the regime of Theorem~\ref{main_thm2_general}, and the isotonic projection of Section~\ref{sec3} is what makes it solvable; we show the worst case collapses to a finite, interpretable \emph{stress ladder} computable in $O(n)$ time.

\textbf{Model.} A platform reserves serving capacity $c\ge 0$ for a future load $X$ with known mean $\mu$ and standard deviation $\sigma$. One capacity unit absorbs $\eta>0$ load units, so the overload after capacity is $X-\eta c$, and a positive overload triggers the contractual service credit. The tiered capped credit is the multi-band utility
\begin{equation}\label{eq:gl-U}
U(y)=\sum_{k=1}^{K} g_k\left (\left (y-a_k\right )_+-\left (y-b_k\right )_+\right ),\quad a_1<b_1\le a_2<b_2\le\cdots\le a_K<b_K,\ g_k>0.
\end{equation}
Band $k$ pays at marginal rate $g_k$ on $(a_k,b_k)$ and then saturates, so the slope of $U$ is the alternating sequence $0,g_1,0,g_2,0,\ldots,g_K,0$. Thus $U$ is nondecreasing and Lipschitz, with a convex kink at each band foot $a_k$ and a concave kink at each band top $b_k$; for $K\ge 2$ it is neither convex nor concave, and it is capped at $\bar U:=U(b_K)=\sum_{k=1}^{K}g_k\left (b_k-a_k\right )$. Write $\bar U_k:=U(b_k)=\sum_{l\le k}g_l\left (b_l-a_l\right )$ for the credit through band $k$, with $\bar U_0:=0$. 

The agreement specifies a finite set of reported percentiles, which we call the named percentiles, so the distortion is purely atomic, placing weight on the quantile levels $0<p_1<p_2<\cdots<p_J<1$, i.e., 
\begin{displaymath}
\mathrm d\varphi(u)=\sum_{j=1}^{J}\omega_j \delta_{p_j}(\mathrm du),\quad \omega_j>0,\ \sum_{j=1}^J\omega_j=1.
\end{displaymath}

For a fixed capacity $c$, the worst case  over all load distributions matching $(\mu,\sigma)$ is
\begin{displaymath}
\Psi_{\max}(c):=\sup_{X\in\mathcal P(\mu,\sigma^2)} \sum_{j=1}^{J}\omega_j U\!\left (F_X^{-1}(p_j)-\eta c\right ).
\end{displaymath}
Here $\mathcal P(\mu,\sigma^2)$ is the set of distributions on $\mathbb R$ with mean $\mu$ and variance $\sigma^2$, and $F_X^{-1}(p):=\inf\{x:F_X(x)>p\}$ is the upper generalized inverse, the convention favorable to the adversary. The worst-case value is the same under the lower inverse, while attainment by a finite staircase uses the upper one. The ambiguity set is taken over $\mathbb R$ rather than the nonnegative axis for tractability, which only enlarges it, so $\Psi_{\max}$ remains a valid and conservative liability bound; because $U$ is nondecreasing and credit accrues only at the upper named percentiles, the negative bulk levels that may appear in the worst case do not affect the credit.
The platform provisions the least capacity for which the worst-case liability does not exceed a service-credit budget $B_0$,
\begin{equation}\label{eq:gl-outer}
\min_{c\ge 0} c, \quad\text{subject to}\ \Psi_{\max}(c)\le B_0.
\end{equation}

Since capacity enters only through the shift $X-\eta c$, the inner value depends on $c$ only through the effective mean $m:=\mu-\eta c$ at fixed variance. Define the reduced worst-case credit
\begin{equation}\label{eq:gl-G}
G(m):=\sup\!\left \{\textstyle\sum_{j=1}^{J}\omega_j\,U(q(p_j)):\ q\in\mathsf K,\ \int_0^1 q(u)\,\mathrm du=m,\ \int_0^1 q^2(u)\,\mathrm du=m^2+\sigma^2\right \}.
\end{equation}
Shifting the feasible set gives the identity $\Psi_{\max}(c)=G(\mu-\eta c)$. Since $G$ is nondecreasing, $\Psi_{\max}$ is nonincreasing in $c$, and $0\le\Psi_{\max}(c)\le\bar U$.

The atomic distortion reduces the infinite-dimensional inner problem to a finite-dimensional one. Define the breakpoint masses
\begin{equation}\label{eq:gl-masses}
\pi_0:=p_1,\quad \pi_j:=p_{j+1}-p_j, \left (1\le j\le J-1\right ),\quad \pi_J:=1-p_J,\quad \textstyle\sum_{j=0}^{J}\pi_j=1.
\end{equation}
For $\sigma^2>0$, Lemma~B.1 in the appendix shows that $G(m)$ equals the $(J{+}1)$-point moment program
\begin{equation}\label{eq:gl-finite}
G(m)=\max\!\left \{\sum_{j=1}^{J}\omega_j U(c_j):\ c_0\le c_1\le\cdots\le c_J,\  \sum_{j=0}^{J}\pi_j c_j=m,\  \sum_{j=0}^{J}\pi_j c_j^2=m^2+\sigma^2\right \},
\end{equation}
in which the support masses $\pi_j$ are fixed by the percentiles and only the levels $c_j$ are free. The maximum is attained, and the optimal levels $(c_0^\ast,c_1^\ast,\ldots,c_J^\ast)$ induce the extremal staircase quantile $q^\ast=c_0^\ast$ on $[0,p_1)$ and $q^\ast=c_j^\ast$ on $[p_j,p_{j+1})$ for $j=1,2,\ldots,J$, with $p_{J+1}:=1$.

Classify each percentile by the credit tier that its optimal level occupies. A percentile $p_j$ is benign if $c_j^\ast<a_1$, so its level earns no credit; breaching if $c_j^\ast\in[a_1,b_K)$, so its level lies in a credit band or on an inter-band plateau; and catastrophic if $c_j^\ast\ge b_K$, so its level attains the cap.

\begin{theorem}[Stress-ladder worst case]\label{thm:gl-ladder}
For $\sigma^2>0$, the worst-case liability equals $\Psi_{\max}(c)=G(\mu-\eta c)$, and it is attained by a staircase load supported on the levels $c_0^\ast\le c_1^\ast\le\cdots\le c_J^\ast$ that solve the finite reduction~\eqref{eq:gl-finite}. The degenerate case $\sigma^2=0$ collapses the ambiguity set to the Dirac mass at $m$, with $G(m)=U(m)$. In every regime with a positive second-moment multiplier, the same levels are obtained as the $\pi$-weighted isotonic regression of $R_c$ over the $J{+}1$ breakpoint levels, computed by the pool-adjacent-violators routine of Algorithm~\ref{algorithm11} with the breakpoint weights $\pi_j$. This worst case has the following structure.
\begin{enumerate}[(a)]
\item There exist multipliers $\lambda_1,\lambda_2\in\mathbb R$, selectors $\xi_j\in\partial U(c_j^\ast)$ in the Clarke subdifferential, and monotonicity multipliers $\nu_j\ge 0$ for $j=0,1,\ldots,J-1$, with the convention $\nu_J:=0$, such that
\begin{equation}\label{eq:gl-kkt0}
0=\pi_0\left (\lambda_1+2\lambda_2 c_0^\ast\right )+\nu_0,
\end{equation}
and
\begin{equation}\label{eq:gl-kkt}
\omega_j\,\xi_j=\pi_j\left (\lambda_1+2\lambda_2 c_j^\ast\right )+\left (\nu_j-\nu_{j-1}\right ),\quad j=1,2,\ldots,J,
\end{equation}
with $\nu_j\left (c_{j+1}^\ast-c_j^\ast\right )=0$ for $j=0,1,\ldots,J-1$. If $\lambda_2>0$, define $\ell^\ast:=-\lambda_1/(2\lambda_2)$; then every percentile whose level is isolated from its neighbors satisfies $c_j^\ast=\ell^\ast+\omega_j\xi_j/(2\lambda_2\pi_j)$, and $c_0^\ast=\ell^\ast$ whenever $c_0^\ast<c_1^\ast$.
\item \textnormal{(Monotone ladder.)} The optimal levels satisfy $c_0^\ast\le c_1^\ast\le\cdots\le c_J^\ast$, so the credit tier of $p_j$ is nondecreasing in $j$. If $c_0^\ast<a_1$, then the bulk mass on $[0,p_1)$ is benign, and the levels $c_j^\ast$ are nondecreasing across the tiers of $U$.
\item \textnormal{(Cap and bound.)} Each catastrophic percentile contributes credit $\bar U$, each benign percentile contributes $0$, and $0\le G(m)\le\bar U$. In positive-multiplier regimes, the projection step is computed in $O(n)$ time by Algorithm~\ref{algorithm11}, as recorded in Remark~\ref{rem:complexity}.
\end{enumerate}
\end{theorem}

By Theorem~\ref{thm:gl-ladder}, the worst-case load is a finite staircase that is benign on its bulk mass and nondecreasing across the named percentiles, with the deepest percentiles on the highest tiers. The levels are not available in closed form when $K\ge 2$, because the monotonicity multipliers $\nu_j$ in Equation \eqref{eq:gl-kkt} couple adjacent percentiles; the characterization instead provides the ordered tier structure and an $O(n)$ projection step in positive-multiplier regimes, without recourse to a convex surrogate. Figure~\ref{fig:gl-stress} visualizes this worst case and makes the input--output distinction explicit: the credit bands of the penalty $U$ appear as the shaded background (the input), while the worst-case load is the staircase $q^\ast$ drawn on top of them (the output), so the band in which each rung lands reads off its benign, breaching, or catastrophic grade.

\begin{figure}[htbp]
\centering
\includegraphics[width=0.92\textwidth]{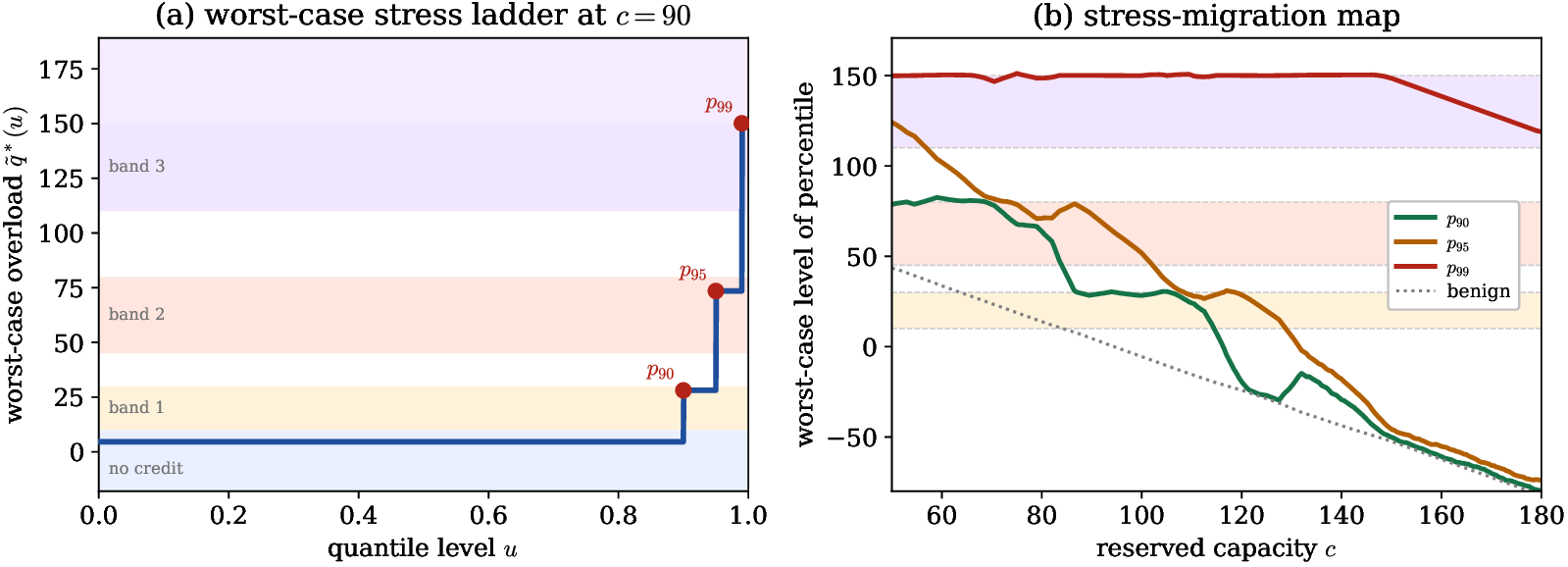}
\caption{The worst-case stress ladder (output) induced by the tiered penalty $U$ and the named-percentile distortion, here at reserved capacity $c=90$. Panel (a): the shaded horizontal bands are the credit tiers of the penalty $U$ (the input---no-credit region, bands $1$--$3$, and the cap), while the blue staircase $q^\ast$ with markers at $p_{90},p_{95},p_{99}$ is the worst-case load (the output); a rung's band is its benign/breaching/catastrophic grade. Panel (b): the stress-migration map---as reserved capacity $c$ falls, the worst-case level of each named percentile climbs through the bands toward the cap.}
\label{fig:gl-stress}
\end{figure}

Since $\Psi_{\max}(c)=G(\mu-\eta c)$ is continuous and nonincreasing in $c$, the provisioning problem~\eqref{eq:gl-outer} reduces to one-dimensional inversion: with $G^{-1}(B_0):=\max\{m\le\mu:G(m)\le B_0\}$, the robust capacity is $c^\star(B_0)=(\mu-G^{-1}(B_0))/\eta$. The price of reliability $|\mathrm dc^\star/\mathrm dB_0|=1/(\eta\,G'(\mu-\eta c^\star))$ reflects the tiered structure of $G$: it is smooth within a tier configuration, jumps when a named percentile crosses a tier, and diverges at a plateau credit level. As reserved capacity falls or volatility grows, the deepest named percentiles climb to the cap. Proposition~B.3 and Corollary~B.1 in the appendix give the frontier, its price, and closed-form levels on each climbing arc.

This tiered frontier is specific to the robust, atomic formulation: a Gaussian fit smooths away the tier structure and under-provisions the deep percentiles, while an absolutely continuous spectral weight of equal mass dissolves the contractual point masses and cannot certify the per-percentile credits the SLA charges on.

\begin{figure}[htbp]
\centering
\includegraphics[width=0.92\textwidth]{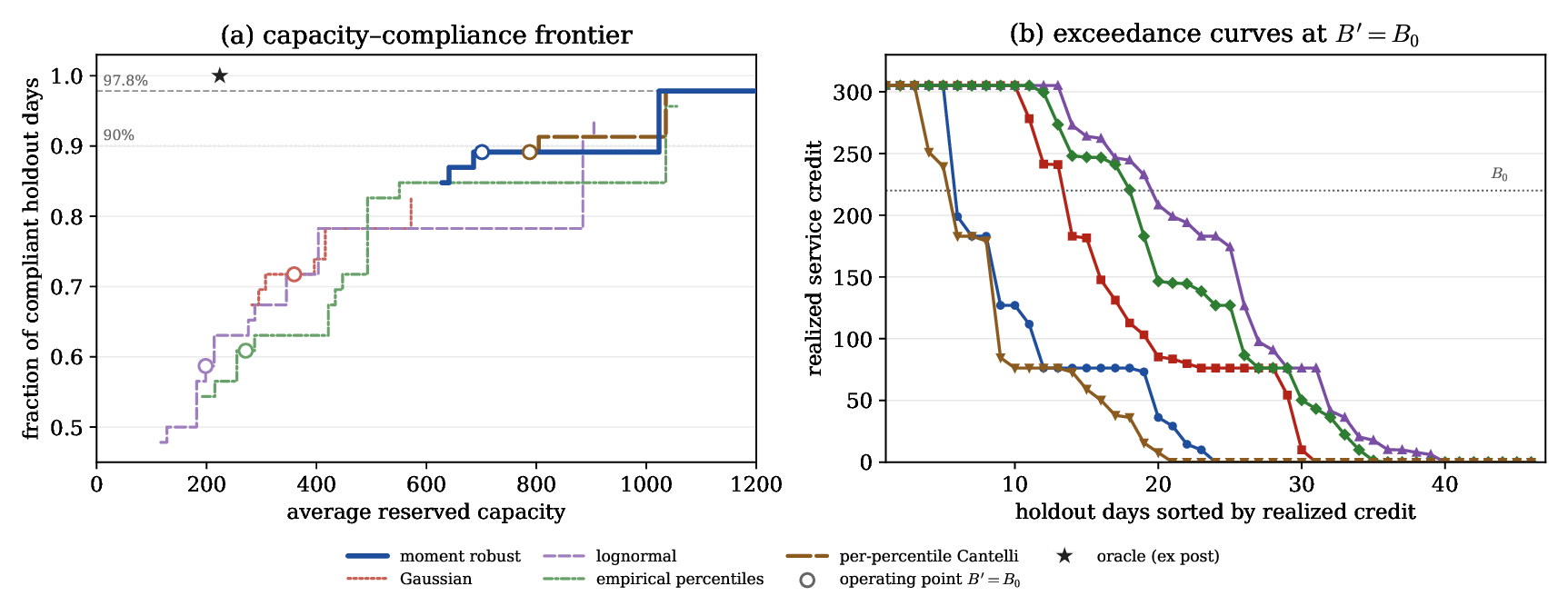}
\caption{Real-trace evaluation on BurstGPT, a public workload trace for large language model serving \citep{wang2025burstgpt}. Load is the five-minute total of request and response tokens, normalized by the full-trace mean. Panel (a) reports capacity-compliance frontiers for budget $B_0=220$: each method varies its internal planning budget, and a holdout day is compliant if realized named-percentile credit is at most $B_0$. Open markers use $B'=B_0$, and the oracle is the ex post minimum daily capacity. The plotted range focuses on average capacities up to 1200; robust and Cantelli also attain full compliance at averages 2313.3 and 2419.4. The dotted and dashed references mark 90.0\% and 97.8\% compliance. Panel (b) sorts realized credits at $B'=B_0$, so budget-line crossings are violation days.}
\label{fig:gl-realtrace}
\end{figure}

Finally, we assess provisioning under estimated moments on a real serving trace. We use 1.40 million successful requests from BurstGPT, aggregate them into five-minute total token loads, and normalize by the full-trace mean. After retaining 60 complete days, the normalized load has mean 101.25, standard deviation 241.59, and 90th, 95th, and 99th percentiles 286.52, 494.65, and 1208.42. This dispersion is where two-moment policies are most exposed to tail misspecification. Each method uses a rolling 14-day window to choose capacity for the next day, leaving 46 holdout days, with $\eta=1$; the internal planning budget $B'$ is the training-window credit ceiling, while out-of-sample compliance is judged against the contractual budget $B_0=220$. Besides the moment-robust plan, we evaluate Gaussian and lognormal fits, training-window empirical percentiles, and a per-percentile Cantelli plan. Appendix~B.11 gives the formal capacity rules and evaluation metrics.

Figure~\ref{fig:gl-realtrace}(a) reports the capacity--compliance frontiers. At the contractual calibration $B'=B_0=220$, the robust plan is compliant on $89.1\%$ of holdout days at average capacity $701.1$, against $71.7\%$ (Gaussian), $58.7\%$ (lognormal), and $60.9\%$ (empirical-percentile). The per-percentile Cantelli plan reaches the same $89.1\%$ on exactly the same five holdout days, but at average capacity $787.9$: the joint characterization of Theorem~\ref{thm:gl-ladder} thus delivers identical realized protection while saving $11.0\%$ of capacity over the textbook two-moment bound that prices each percentile separately. Tightening the planning budget traces the whole frontier---the robust plan reaches $97.8\%$ compliance at average capacity $1023.1$, matching Cantelli while the empirical-percentile plan reaches only $95.7\%$---with the plateau jumps predicted by Proposition~B.3.

The sorted daily credit curves in Figure~\ref{fig:gl-realtrace}(b) show the same mechanism day by day. At $B'=220$ the robust plan leaves only five violation days, against $13$ (Gaussian), $19$ (lognormal), and $18$ (empirical-percentile), and its violation set is contained in each of theirs---a containment forced because those plans price capacity from laws inside the moment ambiguity set, which the robust worst case dominates (Appendix~B.11). The residual violations fall on regime-shift days, when next-day volatility rises well above the training window; the oracle marker, an ex post infeasible benchmark, measures the price of having only training-window moments under this drift. The robust plan also carries the lowest average realized credit, $68.5$ against $118.3$, $159.2$, and $143.9$ for the three alternatives.

\section{Conclusion}\label{sec5}
We studied optimization of DU under moment information, allowing increasing utilities that are nonconvex and nonsmooth together with general distortions. Recasting the problem in the quantile domain, we derived exact first-order optimality conditions and, through isotonic projection of the monotonicity constraint, closed-form extremal values and distributions for both the worst- and best-case problems. The framework is unified and inexpensive: a single algorithm recovers classical moment bounds such as the Scarf bound and extends them to modern DRMs including RVaR, GlueVaR, and the Dual--Power distortion.

Beyond closed forms, the analysis yields decision-relevant structure: the worst-case distribution ranges from low-dimensional discrete laws for piecewise-linear payoffs to continuous distributions for smooth utilities, with the discrete case reducing a robust evaluation to a few interpretable stress scenarios; the extreme-case value is piecewise and monotone in the decision parameters, giving transparent comparative statics; and attainment can fail in tail-heavy regimes. Because the extreme-case DU is exactly the inner problem of a moment-based DRO decision, these characterizations serve as a cheap, reusable inner oracle, illustrated by the capacity-provisioning application where the atomic-distortion characterization lowers required capacity while preserving compliance.

Our analysis characterizes the inner extreme-value problem rather than a full decision problem, and its optimality conditions are local. Natural extensions embed the closed-form bounds into downstream robust decisions, pursue data-driven implementation with estimated moments, and broaden the ambiguity beyond moments to joint distributional and preference ambiguity.

\begingroup
\SingleSpacedXI

\endgroup

\newpage
\pagestyle{fancy}
\fancyhf{}

\lhead{ 
\footnotesize{\textbf{Online Appendices:} \textit{Extreme-Case Distorted Utility under Moment Ambiguity.}  }
}

\rhead{\thepage}
\numberwithin{equation}{section}
\numberwithin{lemma}{section}
\numberwithin{definition}{section}
\numberwithin{proposition}{section}
\clearpage
\pagenumbering{arabic}
\renewcommand*{\thepage}{A-\arabic{page}}
\renewcommand{\thesection}{\Alph{section}}
\renewcommand{\thelemma}{\Alph{section}.\arabic{lemma}}
\renewcommand{\thedefinition}{\Alph{section}.\arabic{definition}}
\renewcommand{\theproposition}{\Alph{section}.\arabic{proposition}}
\renewcommand{\thefigure}{\Alph{section}.\arabic{figure}}
\renewcommand{\thetable}{\Alph{section}.\arabic{table}}
\setcounter{page}{1}
\setcounter{section}{0}
\setcounter{figure}{1}
\setcounter{table}{0}


 \medskip

\end{document}